\documentclass{amsart}

\pdfoutput=1

\usepackage{graphicx}
\usepackage{amsmath}
\usepackage{amsfonts}
\usepackage[caption=false]{subfig}
\usepackage{multirow}

\providecommand{\norm}[1]{\| #1 \|}

\newcommand{\bigpars}[1]{\bigl(#1\bigr)}
\newcommand{\Bigpars}[1]{\Bigl(#1\Bigr)}

\newcommand{\abs}[1]{| #1 |}

\newcommand{\errtable}[4]{
	\begin{tabular}{| c| c | c | c |}
	\hline
	 $\max_i d_i$ & $\tfrac{1}{n}\sum_i d_i$ & $\sqrt{ \tfrac{1}{n} \sum_i d_i^2}$ & \# data points $n$\\
	\hline
	#1 & #2 & #3 & #4 \\
	\hline
	\end{tabular}
	}

\newcommand{\complextable}[3]{
	\begin{tabular}{| c| c | c |}
	\hline
	\# triangles & \# vertices & characteristic length\\
	\hline
	#1 & #2 & #3 \\
	\hline
	\end{tabular}
	}

\title{Simplicial Nonlinear Principal Component Analysis}
\thanks{This work was supported by NSF DMS-1007399}

\author{Thomas Hunt}
\address{Naval Postgraduate School\\ 833 Dyer Road Monterey\\ CA 93943-5216}
\email{twhunt@nps.edu}
\author{Arthur J. Krener}
\address{Naval Postgraduate School\\ 833 Dyer Road Monterey\\ CA 93943-5216}
\email{ajkrener@nps.edu}
\subjclass[2000]{62-04}

\keywords{nonlinear dimensionality reduction, tangent space, manifold learning, principal component analysis}

\begin{document}

\begin{abstract}
We present a new manifold learning algorithm that takes a set of data points lying on or near a lower dimensional manifold as input, possibly with noise, and outputs a simplicial complex that fits the data and the manifold.
We have implemented the algorithm in the case where the input data has arbitrary dimension, but can be triangulated.
We provide triangulations of data sets that fall on the surface of a torus, sphere, swiss roll, and creased sheet embedded in $\mathbb{R}^{50}$.
We also discuss the theoretical justification of our algorithm.
\end{abstract}

\maketitle

%%%%%%%%%%%%%%%%%%%%%%%%%%%%%%%%%%%%%%%%%%%%%%%%%%%%%%%%%%%%%%%%%%%%%%%%%%%%%

\section{Introduction}
Given a large set of data points in a high dimensional space, the task of a manifold learning algorithm is to discover a lower dimensional manifold that approximates the data reasonably well.
Principal component analysis may be interpreted as a manifold learning algorithm when the set of high dimensional data points lie near a lower dimensional affine subspace, in the sense that it can extract the affine subspace from the data set.
Nature does not always serve up inherently linear data sets though.
For example, a set of image vectors generated by photographing a sculpture at different azimuth and altitude angles may intuitively be described by two parameters related to the two angles, but this description is certainly nonlinear \cite{springerlink:10.1007/BF01421486}.
In \cite{MNR:MNR1596}, the authors use principal component analysis to precondition stellar spectra data before passing it to a neural network for classification, and note that a nonlinear preprocessing scheme may help the neural network learn rare or weak features of the data.

Isomap \cite{Tenenbaum22122000}, Local Linear Embedding \cite{Roweis22122000, Saul:2003:TGF:945365.945372}, and Local Tangent Space Alignment \cite{zhang:313} are recent approaches to reducing the dimension of artificially high dimensional data sets without destroying the geometric characteristics of the original data.
Isomap approximates the geodesic distance between every pair of data points by summing the straight line Euclidean distances along the shortest rectilinear path through the data that connects the two points.
It then constructs an embedding of the data into a lower dimensional Euclidean space by applying Multidimensional Scaling to the set of pairwise geodesic distances so that the geodesic distances are nearly preserved in the lower dimensional representation.
Local Linear Embedding approximates each data point as a weighted average of its neighbors, and then reduces the data set by mapping each data point to a lower dimensional space in such a way that the set of nearest neighbors of a data point in the lower dimensional space is preserved under the mapping, and each data point in the lower dimensional space is approximated by the same weighted average of its neighbors.
Given data that lies near the surface of a manifold, Local Tangent Space Alignment estimates the tangent space at each point from its neighbors, and then generates a lower dimensional coordinate system in such a way that the tangent space associated with a point in the lower dimensional space is still aligned with the tangent spaces associated with its neighbors.
It appears that all three of these methods work best when the input data lie on a manifold that admits a global coordinate system into which the data can be mapped.
In the case of Local Linear Embedding, the authors note that is an open question how to modify their algorithm to handle input data that lie on the surface of a sphere or torus \cite[p. 148]{Saul:2003:TGF:945365.945372}.
One of our design motivations was to develop a manifold learning algorithm that can handle this type of data set.

We present Simplicial Nonlinear Principal Component Analysis (SNPCA), an algorithm that shares some of the underlying motivation of Local Tangent Space Alignment, but differs in that it is a manifold learning algorithm whose output is a simplicial complex that acts as a lower dimensional description of the nonlinear input data set.
SNPCA reduces the data set in two senses.
First, every simplex vertex coincides with a surface data point, and typically there will be far fewer simplex vertices than data points.
Secondly, the subset of data points that lie near a face of the simplex are fit by that face, whose dimension is typically much smaller than the dimension of the data points.
We have implemented our algorithm in the case where the data can be fit with a complex of two-simplices, that is, when the data can be triangulated.
This is the case when the data lies near the surface of a two dimensional manifold embedded in $\mathbb{R}^N$.

%%%%%%%%%%%%%%%%%%%%%%%%%%%%%%%%%%%%%%%%%%%%%%%%%%%%%%%%%%%%%%%%%%%%%%%%%%%%%
%%%%%%%%%%%%%%%%%%%%%%%%%%%%%%%%%%%%%%%%%%%%%%%%%%%%%%%%%%%%%%%%%%%%%%%%%%%%%
\section{Algorithm overview for data that can be triangulated}
%%%%%%%%%%%%%%%%%%%%%%%%%%%%%%%%%%%%%%%%%%%%%%%%%%%%%%%%%%%%%%%%%%%%%%%%%%%%%
We initially describe the algorithm for a set of data that can be triangulated, as this is the case for which we have implemented the algorithm. 
The two fundamental inputs into SNPCA are a set of \emph{data vectors} $\{x^i \vert x^i \in \mathbb{R}^N\}$, and a \emph{characteristic length} $\ell$.
The output of SNPCA is a simplicial complex represented by the set $\{ \mathcal{T}, \mathcal{E}, \mathcal{V} \}$, where $\mathcal{V}$ is the set of vertices in the triangulation, and $\mathcal{E}$ and $\mathcal{T}$ are the sets of edges and triangles in the triangulation.
Each triangle in $\mathcal{T}$ is a three element set of vertices from $\mathcal{V}$, and each edge in $\mathcal{E}$ is a two element set of vertices from $\mathcal{V}$.
We assume that the data has been scaled so that changes in the different coordinates of the data vectors are comparably measured.
The algorithm is designed so each vertex coincides with a data point in $\mathbb{R}^N$, but the number of vertices in $\mathcal{V}$ is typically much less than the number of data points.
There are no isolated edges in $\mathcal{E}$, meaning that every edge in $\mathcal{E}$ is a subset of some triangle in $\mathcal{T}$.
There are also no isolated vertices in $\mathcal{V}$, each point in $\mathcal{V}$ is a vertex of at least one edge and one triangle.

The algorithm is composed of two distinct stages.
The first is the \emph{advancing front stage}, where each successful iteration of the main loop yields a new triangle that is added to the triangulation.
In an iteration of the advancing front stage, the algorithm first selects an active \emph{front edge}, which is an edge in $\mathcal{E}$ that belongs to exactly one triangle in $\mathcal{T}$.
A triangle edge can belong to at most two triangles. so intuitively, a front edge is a triangle edge that is exposed.
The algorithm attempts to generate a new triangle composed of the active edge's vertices and another vertex which is typically new, but may be an existing vertex in $\mathcal{E}$.
If the new triangle is acceptable, then the triangulation is updated with the new triangle.
Otherwise, the front edge is unviable for the remainder of the advancing front stage and the algorithm will not return to it until the second stage.
The first stage terminates when every front edge is unviable.

As one can see in the figures in \S \ref{sec:results}, the advancing front stage typically generates an incomplete triangulation of the data.
This incomplete triangulation has \emph{seams}, which are sequences of front edge pairs, where the edges in each pair are both nearly parallel and close to each other.
The second stage of the algorithm is the \emph{seam sewing stage}, where the algorithm attempts to close the remaining gaps in the triangulation.
We now describe both stages in greater detail.

%%%%%%%%%%%%%%%%%%%%%%%%%%%%%%%%%%%%%%%%%%%%%%%%%%%%%%%%%%%%%%%%%%%%%%%%%%%%%
%\subsection{Detailed algorithm description}
%%%%%%%%%%%%%%%%%%%%%%%%%%%%%%%%%%%%%%%%%%%%%%%%%%%%%%%%%%%%%%%%%%%%%%%%%%%%%
\subsection{Advancing front stage for data that is locally two dimensional}
\label{sec:R3_algrthm_scrptn}

%%%%%%%%%%%%%%%%%%%%%%%%%%%%%%%%%%%%%%%%%%%%%%%%%%%%%%%%%%%%%%%%%%%%%%%%%%%%%%%%%%%%%%%%%%%%%%%%%%%%%%%%%%%%%%
\subsubsection{Generating a candidate triangle from the active front edge}
\label{sec:R3_gnrtng_new_trngl}
To generate the candidate triangle based on the active front edge, SNPCA requires the active edge's two vertices, the data points located near each of the two vertices, and the characteristic length.
Once the algorithm has generated the candidate triangle, it then determines whether to add it to the triangulation based on the criteria in \S \ref{sec:accptng_cndt_trngl}.
Given a front edge with vertices $v^1$ and $v^2$, the algorithm determines the best coordinates in $\mathbb{R^N}$ for a third vertex $v^3$ so that the new triangle $\{v^1, v^2, v^3\}$ fits a subset of the data,
where the coordinates are best in the sense that they solve the nonlinear constrained minimization problem (\ref{eqn:cnsrtnd_optmzn_prblm}), which we now describe.

Let $v^1$ and $v^2$ denote the two vertices belonging to the active front edge. 
We define the \emph{empirical local direction covariance matrix} associated with the vertex $v^j$ as
\begin{equation}
\label{eqn:emprcl_drctn_mtrx}
P(v^j, \mathcal{N}(v^j))
\equiv
\frac{1}{\abs{\mathcal{N}(v^j)}}\sum_{x^i \in \mathcal{N}(v^j) } \frac{(x^i-v^j)(x^i-v^j)^T}{(x^i - v^j)^T(x^i - v^j)} 
\end{equation}
where $v^j$ and $x^i$ are column vectors of coordinates, $\mathcal{N}(v^j)$ is a subset of data falling in a neighborhood centered at the vertex $v^j$, and $\abs{\mathcal{N}(v^j)}$ denotes the number of data points in the neighborhood.
Two possibilities for $\mathcal{N}(v^j)$ are a Euclidean neighborhood, or the $k$ data points nearest $v^j$.

The definition of the empirical local direction matrix is motivated by ideas from principal component analysis.
If the vertex $v^j$ and the data points  in $\mathcal{N}(v^j)$ lie on the surface of a smooth two dimensional manifold, then each member of the set of normalized directions $\{(x^i - v^j)/\norm{x^i-v^j}_2 \mid x^i \in \mathcal{N}(v^i) \}$
can be nearly reconstructed from a two dimensional basis, although the set of directions almost certainly spans a much higher dimensional space.
Akin to principal component analysis, the dominant eigenvectors of the covariance matrix $P(v^j, \mathcal{N}(v^j))$ furnish a natural basis for the set of normalized directions, and the eigenvalues provide a measure of how well each basis vector reconstructs the data.

Associated with $P(v^j, \mathcal{N}(v^j))$ is a Riemannian metric induced by the the matrix
\begin{equation}
\label{eqn:emprcl_drctn_mtrx_invrs}
Q_{\mu}(v^j, \mathcal{N}(v^j)) \equiv \bigpars{P(v^j, \mathcal{N}(v^j)) + \mu I}^{-1}
\text{, }  \mu \ge 0
\end{equation}
where the user must set the parameter $\mu$ so that under this metric, the distance from $x$ to $v^j$ is large when $x$ does not lie in the affine subspace associated with dominant eigendirections of $P(v^j, \mathcal{N}(v^j))$.
Without an efficient method to compute the distance under the metric induced by $Q_{\mu}(v^j, \mathcal{N}(v^j))$, our algorithm would require an impractical amount of time to complete.
Our assumption that the data can be fit with triangles implies that the numerical rank of $P_{\mu}(v^j, \mathcal{N}(v^j))$ is much less than the dimension of a data vector, so we only need the dominant eigenvalues and eigenvectors of $P(v^j, \mathcal{N}(v^j))$ to measure distances in the induced metric.
Suppose $\mu >0$ and $P_{\mu}(v^j, \mathcal{N}(v^j))$ has $k<N$ nonzero eigenvalues. 
Let $\Lambda$ denote the $k \times k$ diagonal matrix of nonzero eigenvalues, and let $V$ denote the $N \times k$ matrix whose orthonormal columns are eigenvectors of $P_{\mu}(v^j, \mathcal{N}(v^j))$ so that $P_{\mu}(v^j, \mathcal{N}(v^j))V = V \Lambda$.
Let $y \in \mathbb{R}^k$ be the coordinate vector $y \equiv V^Tx$ in, then the distance induced by $Q_{\mu}(v^j, \mathcal{N}(v^j))$ \eqref{eqn:emprcl_drctn_mtrx_invrs} can be computed as
\begin{equation} \label{eqn:Q_inner_prod_fast}
x^T\bigpars{P_{\mu}(v^j, \mathcal{N}(v^j)) + \mu I}^{-1}x
=
y^T\bigpars{\Lambda + \mu I}^{-1} y + \frac{1}{\mu}\bigpars{\norm{x}_2^2 - \norm{y}_2^2}
\end{equation}
Computing the lefthand side of \eqref{eqn:Q_inner_prod_fast} in the naive but straightforward manner requires $\mathcal{O}(N^2)$ arithmetic operations.
Forming $y$ requires $\mathcal{O}(k N)$ arithmetic operations, and once we have formed $y$, forming the righthand side of \eqref{eqn:Q_inner_prod_fast} requires an additional $\mathcal{O}(N)$ arithmetic operations.
Of course, the righthand side requires the eigen information in $V$ and $\Lambda$, but these can be computed efficiently by a Krylov subspace method since the structure of $P_{\mu}(v^j, \mathcal{N}(v^j))$ allows us to compute the matrix-vector product $P_{\mu}(v^j, \mathcal{N}(v^j)) z$ as the linear combination of the normalized directions $x^i-v^j$.
This requires $\mathcal{O}(\abs{\mathcal{N}(v^j)} N)$ arithmetic operations, where $\abs{\mathcal{N}(v^j)}$ is the number of normalized directions, and typically $\abs{\mathcal{N}(v^j)}$ is much less than the dimension of a data vector $N$.
 
Given the front edge with vertices $v^1$ and $v^2$, the algorithm attempts to generate a new triangle with the vertices $\{v^1, v^2, v^3 \}$, where $v^3$ is typically new, but may be an existing vertex.
The bulk of the computational expense required to generate the new triangle is spent solving the constrained minimization problem (\ref{eqn:cnsrtnd_optmzn_prblm}) for the new vertex.
Let $Q^j$ denote $Q(v^j, \mathcal{N}(v^j))$, then the constrained optimization problem is
\begin{subequations}
\label{eqn:cnsrtnd_optmzn_prblm}
\begin{align}
v^*
&= 
\underset{\mathcal{S}}
{\text{argmin}}  
(v-v^1)^T Q^1 (v-v^1)
\label{eqn:opmtzn_objctv}
\\
\begin{split}
\label{eqn:opmtzn_cnstrnts}
\mathcal{S}
&:
\text{$v$ in $\mathbb{R}^3$ such that}
\left\{
\begin{aligned}
(v - v^m)^T(v - v^m) &= R_c^2 \\
(v-v^1)^T Q^1 (v-v^1) &= (v-v^2)^T Q^2 (v-v^2) \\
n^T(v-v^m) &\ge 0 
\end{aligned}
\right.
\end{split}
\end{align}
\end{subequations}

The first constraint describes a Euclidean sphere of radius $R_c$ centered at the midpoint $v^m \equiv (v^1 + v^2)/2$ of the active front edge.
We assume that the raw data has been scaled so that all the points on this sphere represent more or less equal  
changes in the characteristics of the data.
The second constraint is the isosceles constraint, so called because it forces the minimizer to be equidistant from both active edge vertices, where the distance to the edge vertex $v^j$ is measured by the metric induced by $Q^j$.
The shorter the $Q^j$ length of a displacement  from $v^j$, the more confident we are that the data lies in that direction.
Under the first two constraints in (\ref{eqn:opmtzn_cnstrnts}), the minimization problem can have two solutions.
Setting the new vertex to one of these minimizers would introduce a triangle that is redundant in the sense that it nearly replicates the triangle that owns the active edge.
The third constraint is designed to remove this minimizer from consideration.
Let $v^0$ denote the third vertex in the triangle that owns the active front edge, then $n$ is the unit normal to the front edge that lies in the plane of $\{v^0, v^1, v^2\}$, and points away from $v^0$.
That is, it is the unit vector in the direction of $(v^m-v^0) - ((v^m-v^0)^T(v^2-v^1)/\norm{v^2-v^1}^2)(v^2-v^1)$.

\begin{figure}
\subfloat[First candidate vertex]{\includegraphics[scale=.4]{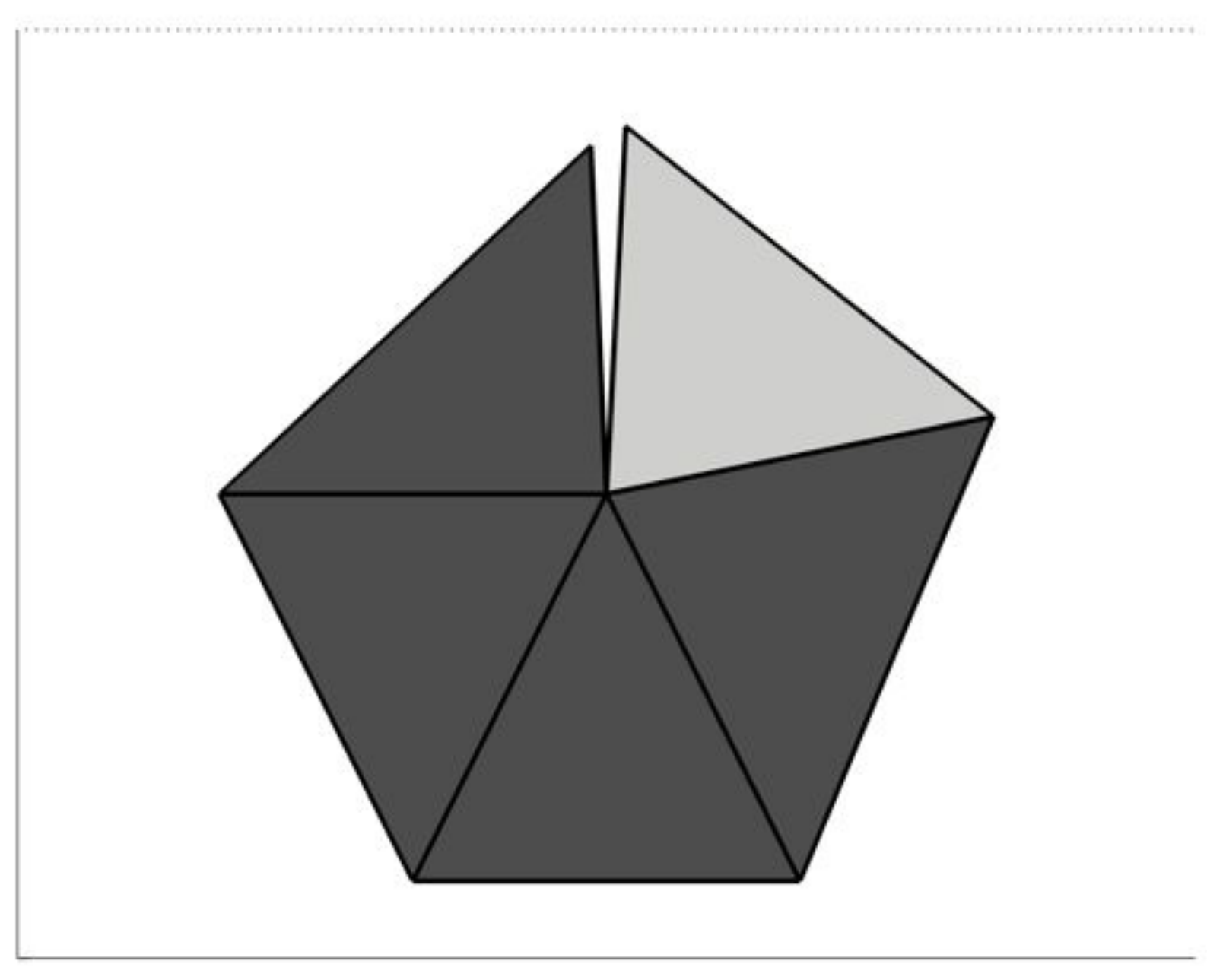}}
\subfloat[Second candidate vertex]{\includegraphics[scale=.4]{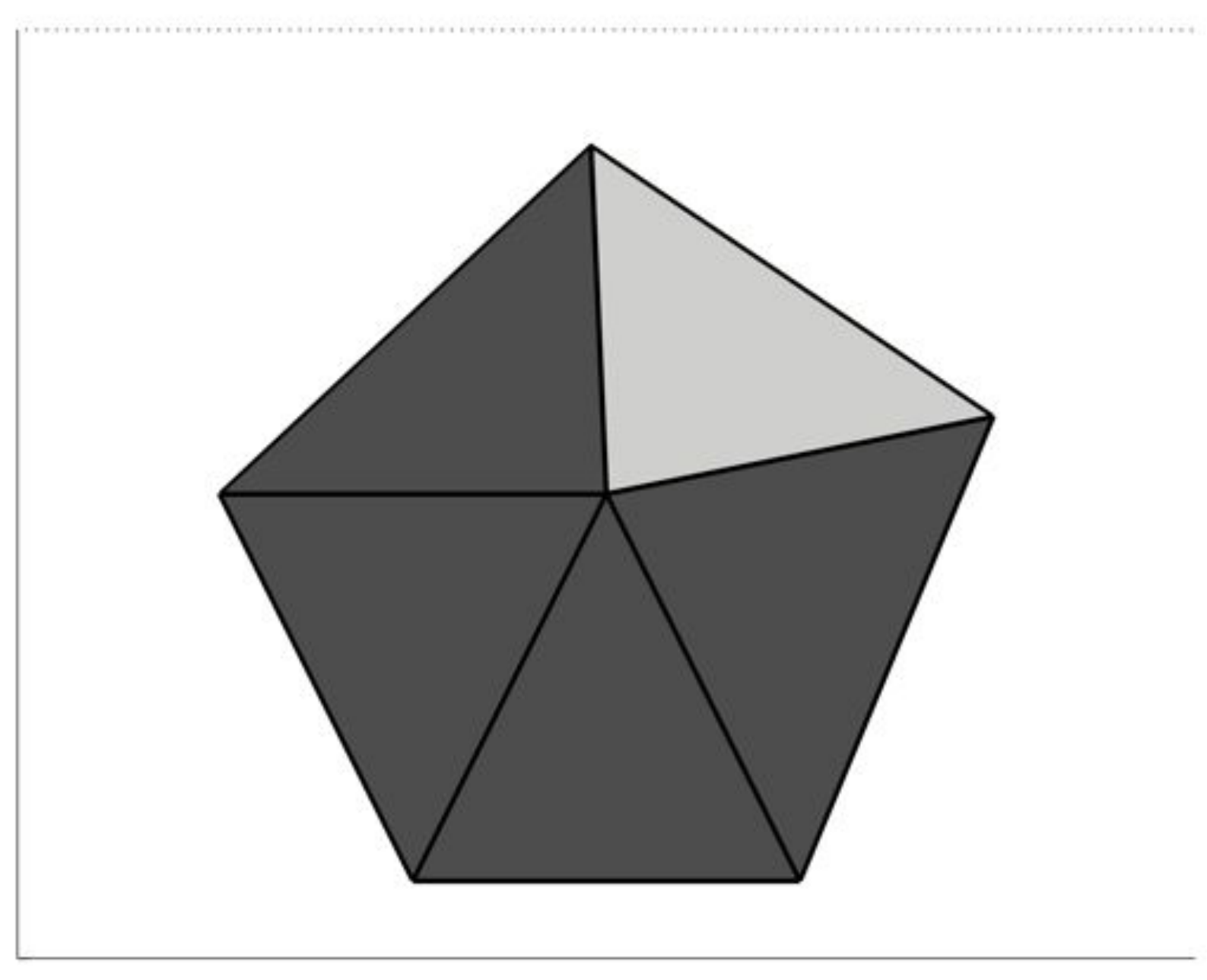}}
\caption{First and second candidate candidate vertices. Existing triangles are dark gray, candidate triangle is light gray}
\label{fig:cndt_vrtcs}
\end{figure}

After solving the constrained optimization problem (\ref{eqn:cnsrtnd_optmzn_prblm}), the algorithm uses the minimizer $v^*$ to generate at most two candidates for the vertex $v^3$.
The first candidate always exists, and is the data point nearest $v^*$, which may also be a vertex in the existing triangulation.
The purpose of the second candidate is to avoid generating an unnecessarily small triangle at some point in the future, which works against the goal of fitting the data with as few triangles as possible.
This can happen if the first candidate vertex is close to a vertex belonging to some other front edge that shares a vertex with the active front edge.
In this event, the second candidate for $v^3$ is the nearby vertex belonging to the front edge adjacent to the active edge, provided that the distance between the two vertices is within a user supplied tolerance.
Without considering the second candidate, then the original candidate triangle will have a front edge that forms a small angle with an adjacent front edge.
As a result, the algorithm will likely generate a future small triangle that owns the two adjacent front edges with a small angle between them, and a third short edge.
Illustrations of the first and second candidate vertices are in Figure \ref{fig:cndt_vrtcs} in the case where the data vectors are in $\mathbb{R}^3$.

The user must supply the radius of the constraint sphere $R_c$ in (\ref{eqn:opmtzn_cnstrnts}).
In our implementation, we set the constraint sphere radius based on the weighted average of the Euclidean length of the active front edge, and the characteristic length.
A weighted average with a heavy emphasis on the active edge length is more likely to result in an equilateral triangle, which is preferable if the active front edge is not too short.
If the active front edge is relatively short, then it is preferable to generate an isosceles triangle with two long edges and a single short edge since this triangle will fit more data than an equilateral triangle.
In our implementation, we set $R_c=\sqrt{3}/2$ times the weighted average of the active edge length and the characteristic length, where both weights are $1/2$.
In the common case where the active front edge length and the characteristic length are nearly equal, then this constraint sphere radius will favor a Euclidean equilateral triangle with sides that are approximately the characteristic length.

In our implementation, we use a Euclidean neighborhood for $\mathcal{N}(v^i)$.
When calculating $P(v^i, \mathcal{N}(v^i))$, we set the radius of this neighborhood to the length of the active edge $\{v^1, v^2\}$ measured in the Euclidean metric.

Given the popularity and long history of least squares, one may ask why we did not choose to minimize the sum of squares $(v-v^1)^T Q^1 (v-v^1) + (v-v^2)^T Q^2 (v-v^2)$ instead of the objective function \eqref{eqn:opmtzn_objctv} in the constrained minimization problem.
It turns out that when $Q^1 = Q^2 = Q$, as may be the case when the data is planar, the least squares solution is unacceptable.
Suppose the origin is the midpoint of the active front edge, so its vertices satisfy $v^2 = -v^1$.
Then, $(v-v^1)^T Q (v-v^1) + (v-v^2)^T Q (v-v^2) = 2 x^T Q x + 2 x^T Q x^1$, which is minimized on the constraint sphere \eqref{eqn:opmtzn_cnstrnts} by the appropriately scaled eigenvector of $Q$ associated with its smallest eigenvalue.
In this case, the least squares solution completely ignores the other two eigenvalues of $Q$, so this information is essentially wasted.
Furthermore, there is no reason to believe that the minimizing eigenvector of $Q$ is nearly perpendicular to the active front edge, so the resulting triangle may be highly skewed.

%%%%%%%%%%%%%%%%%%%%%%%%%%%%%%%%%%%%%%%%%%%%%%%%%%%%%%%%%%%%%%%%%%%%%%%%%%%%%%%%%%%%%%%%%%%%%%%%%%%%%%%%%%%%
\subsubsection{Criteria for accepting the candidate triangle}
\label{sec:accptng_cndt_trngl}
The main idea behind deciding whether to accept the candidate triangle is that the new triangle should not \emph{conflict} with the existing triangulation, which we now define.
We can always translate the triangles $T^1$ and $T^2$ so that without a loss of generality, we may assume that one of the vertices of $T^1$ coincides with the origin in $\mathbb{R}^N$.
Let $\bar{T}^1$ and $\bar{T}^2$ denote the triangles in $\mathbb{R}^2$ obtained by orthogonally projecting  the vertices of $T^1$ and $T^2$ into the subspace containing the vertices of $T^1$.
Note that $\bar{T}^1$ and $T^1$ represent the same triangle, but their vertices are represented in different bases, while $\bar{T}^2$ and $T^2$ do not represent the same triangle in general.
The triangles $T^1$ and $T^2$ are said to \emph{overlap} if the intersection of the open interiors of $\bar{T}^1$ and $\bar{T}^2$ is nonempty.
More concretely, let $\bar{\mathcal{H}}^1$ and $\bar{\mathcal{H}}^2$ denote the open interiors of the triangles $\bar{T}^1$ and $\bar{T}^2$, that is $\bar{\mathcal{H}}^2 \equiv \{ t_1 (\bar{w}^1-\bar{w}^0) + t_2(\bar{w}^2-\bar{w}^0) \mid t_1, t_2 >0; 0 < t_1 + t_2 < 1 \}$,  $\bar{w}^i$ is a vertex of $\bar{T}^2$, and $\bar{\mathcal{H}}^1$ is defined analogously.
$T^1$ and $T^2$ overlap if $\bar{\mathcal{H}}^1 \cap \bar{\mathcal{H}}^2$ is nonempty.
The triangles $T^1$ and $T^2$ are said to \emph{conflict} if they overlap and the distance between $T^1$ and $T^2$ is within some tolerance.
We included the open instead of closed interior in the definition of overlap so that two triangles that share a common edge do not necessarily conflict.

\begin{figure}
\begin{centering}
\subfloat[Overlapping but nonconflicting triangles]
    {\includegraphics[scale=.5]{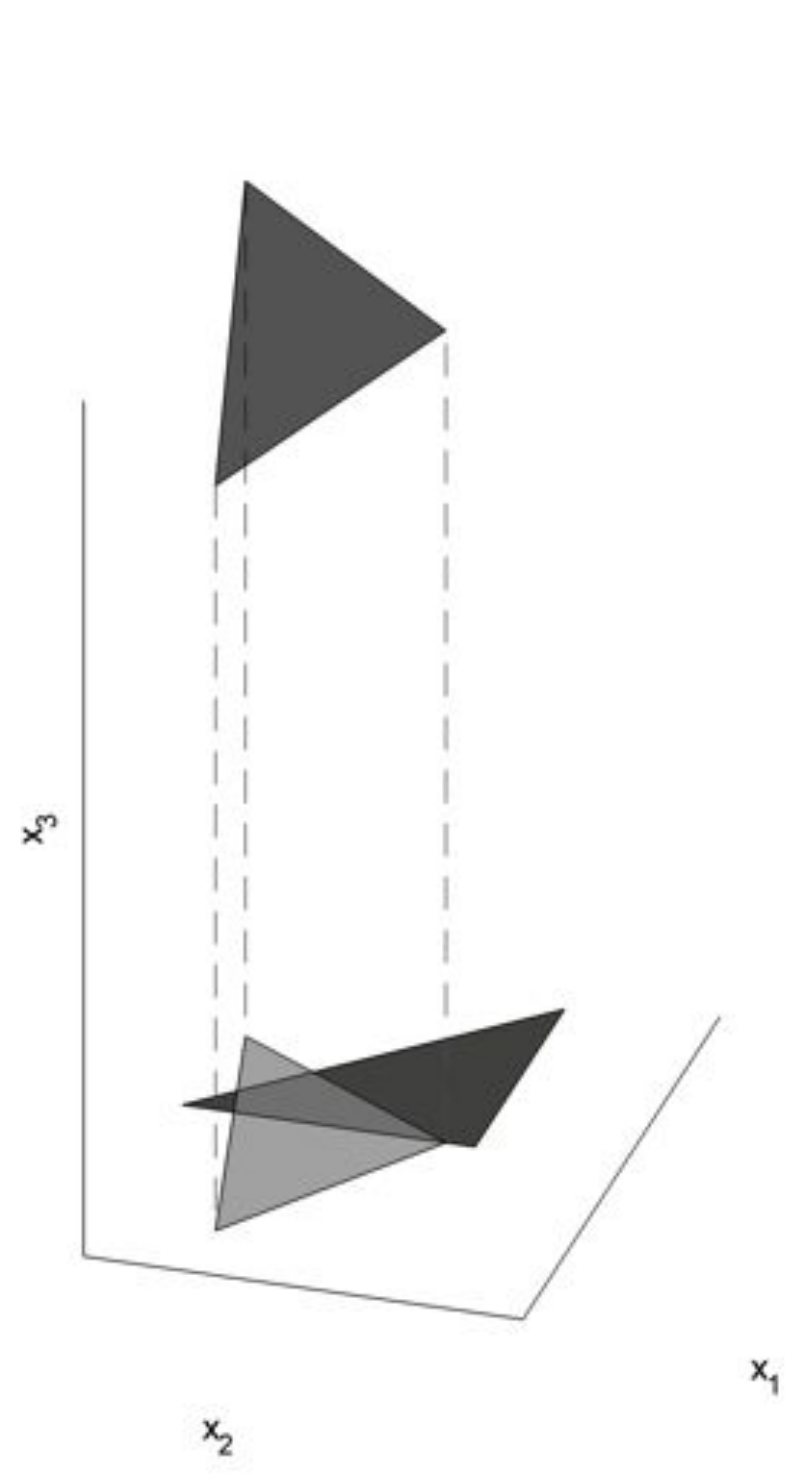}}
\subfloat[Conflicting triangles]{\includegraphics[scale=.5]{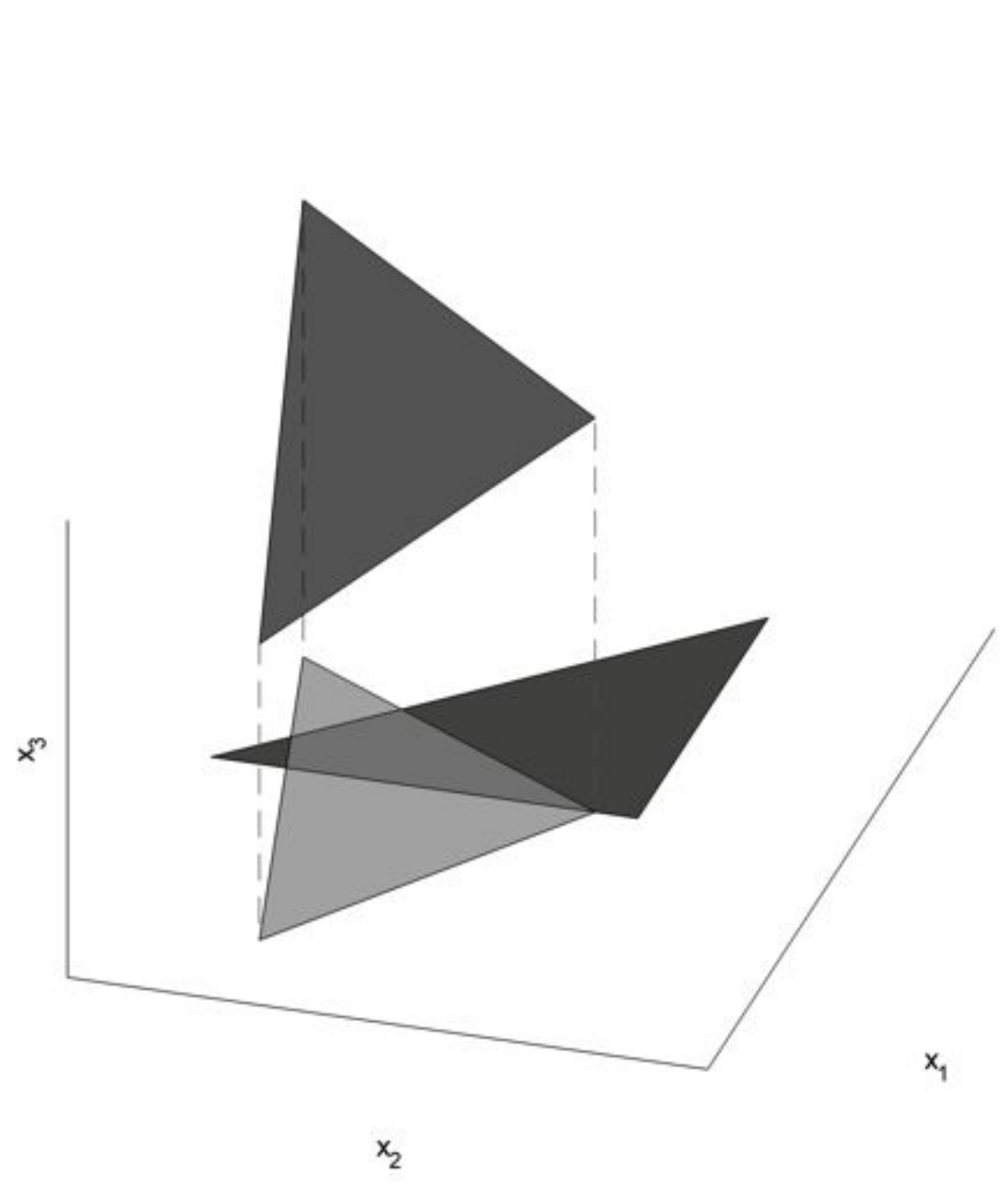}}
\end{centering}
\caption{Conflicting and nonconflicting triangles}
\label{fig:tri_tri_cnflct}
\end{figure}

Figure \ref{fig:tri_tri_cnflct} illustrates the two cases where a pair of triangles overlap but do not conflict, and when a pair of triangles do conflict.
The illustrations are for data in $\mathbb{R}^3$ for the purpose of visualization, but the definitions of conflict and overlap hold for data in a higher dimensional space.
In both cases, triangle $T^2$ is the dark triangle highest on the vertical axis, $T^1$ is the dark triangle in the $x_1$-$x_2$ plane, and the lighter triangle is the triangle formed by orthogonally projecting the vertices of $T^2$ into the plane containing $T^1$.

If the candidate triangle conflicts with an existing triangle in the triangulation, then the algorithm rejects the candidate triangle and adds no new vertices or edges to the triangulation in the current iteration of the main loop.
The algorithm may also fail to generate a new triangle if
\begin{itemize}
\item
The matrix $Q^j$ could not be computed for some active front edge vertex $v^j$ because the deleted neighborhood of $v^j$ contained no surface data points.
This indicates that the radius of the neighborhood $\mathcal{N}(v^j)$ may be too small, or the density of data points near $v^j$ is too low.

\item
The distance from the minimizer $v^*$ to the nearest surface data point was larger than a user supplied tolerance.
This can happen when a front edge lies near a boundary of the data set, such as the boundary of the creased sheet.

\item
The constraint set $\mathcal{S}$ in (\ref{eqn:opmtzn_cnstrnts}) was empty. 
\end{itemize}

%%%%%%%%%%%%%%%%%%%%%%%%%%%%%%%%%%%%%%%%%%%%%%%%%%%%%%%%%%%%%%%%%%%%%%%%%%%%%
\subsubsection{Generating the initial triangulation}
\label{sec:intl_prtl_trngltn}
The advancing front stage of the algorithm requires an initial triangulation to get started.
The algorithm generates the first vertex $v^1$ by setting its coordinates to those of some data point that is chosen randomly unless the user specifies a specific initial point.
The algorithm then finds a second data point so that the distance from the initial vertex $v^1$ to this data point is as close to the characteristic length as possible, and sets the coordinates of the second vertex $v^2$ equal to the coordinates of the second data point.
With the edge defined by vertices $v^1$ and $v^2$, the algorithm now solves the constrained minimization problem of (\ref{eqn:cnsrtnd_optmzn_prblm}) with the third constraint removed, which typically has two solutions.
The algorithm generates coordinates for two new vertices from these minimizers by finding the two data points closest to the minimizers, and then generates the initial triangulation containing four vertices, five edges, and two triangles.

%%%%%%%%%%%%%%%%%%%%%%%%%%%%%%%%%%%%%%%%%%%%%%%%%%%%%%%%%%%%%%%%%%%%%%%%%%%%%
\subsubsection{Selecting the active front edge in the advancing front stage}
\label{sec:chsng_actv_frnt_edg}
The algorithm must select the active front edge at the beginning of each iteration of the advancing front stage of the algorithm.
In our implementation, we push any new front edges onto a stack (FILO buffer) whenever a new triangle is created.
Although an edge may be a front edge when it is pushed onto the stack, it may get covered by a triangle after some subsequent iteration of the advancing front stage.
So, at the beginning of the advancing front stage, the algorithm pops edges off the stack until it encounters a front edge, which becomes the active front edge.
This approach appears to work fairly well in practice, although a more sophisticated active front edge choosing algorithm may yield a better final triangulation.

%%%%%%%%%%%%%%%%%%%%%%%%%%%%%%%%%%%%%%%%%%%%%%%%%%%%%%%%%%%%%%%%%%%%%%%%%%%%%
\subsection{Seam sewing stage for data that can be triangulated}
\label{sec:seam_sewing_stage}
The advancing front stage of SNPCA typically produces a triangulation that does not completely fit the data set, as one can see in \S \ref{sec:results}.
The seam sewing stage of the algorithm attempts to complete the triangulation with a variation of the advancing front stage algorithm.
This stage requires the partial complex $\{ \mathcal{T}, \mathcal{E}, \mathcal{V} \}$ generated during the advancing front stage, the set of $Q^j$'s associated with all vertices belonging to a front edge, and a maximum allowable edge length parameter that is similar in spirit to the characteristic length.
Each successful iteration of the seam sewing stage produces a new triangle from an active edge,
but unlike the advancing front stage, it never generates a new vertex.

%%%%%%%%%%%%%%%%%%%%%%%%%%%%%%%%%%%%%%%%%%%%%%%%%%%%%%%%%%%%%%%%%%%%%%%%%%%%%
\subsubsection{Generating a new triangle from the active front edge}
Let $v^1$ and $v^2$ denote the vertices of the active front edge.
The algorithm attempts to form a new triangle $\{v^1, v^2, v^3\}$ by finding the existing vertex $v^3$ that solves the discrete optimization problem
\begin{equation}
\label{eqn:sew_seams_optmzn_prblm}
\begin{split}
v^3
&= 
\underset{\mathcal{V}}
{\text{argmin}}  
(v-v^1)^T Q^1 (v-v^1) + (v-v^2)^T Q^2 (v-v^2)
\\
&\text{where}
\left\{
\begin{aligned}
&\norm{v^3 - v^1}_2, \norm{v^3 - v^2}_2 \le \text{ maximum allowable length} \\
&\{v^1, v^2, v^3\} \text{ does not conflict with the existing triangulation}
\end{aligned}
\right.
\end{split}
\end{equation}
Note that the minimizer $v^3$ comes from the finite set of existing vertices $\mathcal{V}$, so the first constraint reduces the problem to a search over a small subset of the vertex set.
Since this set is discrete, there is no reason to believe that the isosceles constraint of (\ref{eqn:opmtzn_cnstrnts}) can be satisfied.
Consequently, we chose the objective function to simply minimize the sum of the squares of the triangle edge lengths opposite the active edge, measured in their respective induced metrics.
We imposed the first constraint so the seam sewing stage of the algorithm does not generate a new triangle that bridges a void in the data.
The algorithm marks the active front edge as unviable if no vertex satisfies the constraints, and this information is used by the algorithm when selecting the next active edge in subsequent iterations.

In practice, the seam sewing stage may fail to close holes in the triangulation, in which case it it is up to the user to determine if the holes are spurious.
If so, they may be closed by omitting the restriction in (\ref{eqn:sew_seams_optmzn_prblm}) that the candidate triangle does not conflict with the existing triangulation.

%%%%%%%%%%%%%%%%%%%%%%%%%%%%%%%%%%%%%%%%%%%%%%%%%%%%%%%%%%%%%%%%%%%%%%%%%%%%%
\subsubsection{Selecting the active edge in the seam sewing stage}
\label{sec:seam_sewing_chsng_actv_frnt_edg}
The algorithm that selects the next active edge during the seam sewing stage of the algorithm is similar to the one used during the advancing front stage.
The major difference being that the advancing front stage produced a single sequence of front edges that terminated as soon as the edge stack became empty, whereas the seam sewing stage may produce more than one such sequence. 
At the beginning of the seam sewing stage, each front edge is marked as viable, and a front edge is marked unviable if upon visiting it, the algorithm determines that the minimization problem (\ref{eqn:sew_seams_optmzn_prblm}) has no solution.

To initialize a sequence of front edges, the algorithm finds the smallest angle between two adjacent front edges, where at least one of the front edges is viable, and then pushes one of the viable front edges into the edge stack.
The algorithm repeatedly pops edges out of the buffer until it finds a front edge that has not been marked as unviable, and then sets the active front edge to this edge.
It then attempts to generate a new triangle by solving the minimization problem in (\ref{eqn:sew_seams_optmzn_prblm}), and adds any new front edges to the edge stack.
The sequence of front edges terminates when the edge stack is empty, and the seam sewing stage terminates when there are no front edges in the triangulation or every front edge is unviable.

%%%%%%%%%%%%%%%%%%%%%%%%%%%%%%%%%%%%%%%%%%%%%%%%%%%%%%%%%%%%%%%%%%%%%%%%%%%%%
\subsection{Advancing front stage outline for data that can be fit with a complex of $d$-simplices}
In this section, we give an overview of the general case where the data points are in $\mathbb{R}^N$, but are essentially locally $d$-dimensional.
This means that given any data point, the data points that fall in a sufficiently small neighborhood of that data point lie near a $d$-dimensional affine subspace.
In this setting, the output of SNPCA is a simplicial $d$-complex whose simplices have at most $d+1$ vertices.
If the user does not know the essential dimension of the data beforehand, it can be estimated at a data point by counting the number of dominant eigenvalues of the empirical local direction covariance matrix computed at that data point.

Analogous to a front edge, a \emph{front face} is a $d-1$-simplex in the complex whose vertices are a subset of exactly one $d$-simplex in the complex.
The fundamental task in the advancing front stage is to generate a new vertex from a front face, so that the new vertex and front face form a new $d$-simplex that fits a subset of the data, which is a generalization of generating a new triangle from a front edge in \S \ref{sec:R3_gnrtng_new_trngl}.
Let $\{v^1, v^2, \ldots, v^d\}$ denote the \emph{active front face} vertices that belong to the $d$-simplex with vertices $\{v^0, v^1, v^2, \ldots, v^d\}$, define the local direction covariance matrix  and $Q_{\mu}(v^j, \mathcal{N}(v^j))$ as in (\ref{eqn:emprcl_drctn_mtrx}) and (\ref{eqn:emprcl_drctn_mtrx_invrs}).
The coordinates for the new vertex are generated by first solving the following generalization of the constrained minimization problem (\ref{eqn:cnsrtnd_optmzn_prblm})
\begin{subequations}
\label{eqn:gnrl_cnsrtnd_optmzn_prblm}
\begin{align}
v^*
&= 
\underset{\mathcal{S}}
{\text{argmin}}  
(v-v^1)^T Q^1 (v-v^1)
\label{eqn:gnrl_opmtzn_objctv}
\\
\mathcal{S}
&:
\text{$v$ in $\mathbb{R}^N$ such that}
\left\{
\begin{aligned}
(v - v^m)^T(v - v^m) &= R_c^2 \\
(v-v^1)^T Q^1 (v-v^1) &= (v-v^j)^T Q^j (v-v^j), \; 2 \le j \le n \\
n^T(v-v^m) &\ge 0 
\end{aligned}
\right.
\label{eqn:gnrl_opmtzn_cnstrnts}
%\end{split}
\end{align}
\end{subequations}
where $v^m \equiv (v^1 + v^2 + \cdots + v^d)/d$ denotes the centroid the active front face, $Q^j$ is shorthand for $Q_{\mu}(v^j, \mathcal{N}(v^j))$, and $R_c$ denotes the radius of the constraint sphere.
The second constraint forces the minimizer $v^*$ to be equidistant from all front face vertices where the distance to the $j^{\text{th}}$ vertex is measured with the metric induced by $Q^j$.
The purpose of the third constraint in (\ref{eqn:gnrl_opmtzn_cnstrnts}) is to force the minimizer to fall on the side of the active face opposite the interior vertex $v^0$, which is the intuitively correct region to place the vertex.
We take $n$ to be the unit normal to the $(d-1)$-simplex that is the active front face, lies in the space spanned by the simplex that owns the front face, and points away from $v^0$.

The other main component of the advancing front stage is a means to determine if a candidate $d$-simplex conflicts with an existing \mbox{$d$-simplex} in the complex, which can be naturally generalized from the definition of conflicting triangles in \S \ref{sec:accptng_cndt_trngl}.
Let $\mathcal{S}^1$ and $\mathcal{S}^2$ denote two simplices whose vertices are in $\mathbb{R}^N$, where as in the case of triangles we may assume that $\mathcal{S}^1$ has one vertex at the origin.
Let $\bar{\mathcal{S}}^1$ and $\bar{\mathcal{S}}^2$ denote the $d$-simplices in $\mathbb{R}^d$ obtained by orthogonally projecting the vertices of $\mathcal{S}^1$ and $\mathcal{S}^2$ into the subspace containing the vertices of $\mathcal{S}^1$.
The simplices $\mathcal{S}^1$ and $\mathcal{S}^2$ are said to overlap if the intersection of the open interiors of $\bar{\mathcal{S}}^1$ and $\bar{\mathcal{S}}^2$ is nonempty.
$\mathcal{S}^1$ and $\mathcal{S}^2$ are said to conflict if they overlap and the distance from $\mathcal{S}^1$ to $\mathcal{S}^2$ is within some tolerance.

If the candidate $d$-simplex formed by the active front face and the new vertex does not conflict with an existing $d$ simplex in the complex, then the algorithm accepts the new simplex.
Just as in the case where the data could be triangulated, the algorithm pushes all new front faces belonging to the new $d$-simplex onto a face stack, and at the top of each iteration of the advancing front stage main loop, faces are popped off the stack until the algorithm finds a front face.

%%%%%%%%%%%%%%%%%%%%%%%%%%%%%%%%%%%%%%%%%%%%%%%%%%%%%%%%%%%%%%%%%%%%%%%%%%%%%
\section{Discussion}
%%%%%%%%%%%%%%%%%%%%%%%%%%%%%%%%%%%%%%%%%%%%%%%%%%%%%%%%%%%%%%%%%%%%%%%%%%%%%
%%%%%%%%%%%%%%%%%%%%%%%%%%%%%%%%%%%%%%%%%%%%%%%%%%%%%%%%%%%%%%%%%%%%%%%%%%%%%%%%%%%%%%%%%%%%%%%%%%%%%%%%%%
\subsection{Interpreting the local direction covariance matrix}
\label{sec:lcl_cvrnc_mtrx_intrptn}
Although the SNPCA algorithm operates on discrete data, we can gain intuition about its behavior by investigating some of its components in a continuous setting.
To this end, suppose we have a smooth manifold $x:\mathbb{R}^2 \rightarrow \mathbb{R}^3$ that can be parameterized locally by the variables $s$ and $t$.
This is the limiting case where the data is distributed uniformly over the manifold with respect to surface area, and the number of data points goes to infinity.
We define the continuous version of the local direction covariance matrix as
\begin{equation}
\label{eqn:exct_emprcl_drctn_mtrx}
P(\bar{v}, \mathcal{N}) \equiv \int_{\mathcal{N}} \frac{(x - \bar{v})(x - \bar{v})^T}{(x - \bar{v})^T(x - \bar{v})} dS
\Bigg/
\int_{\mathcal{N}} dS
\end{equation}
where $\bar{v}$ is a point on the manifold, and $\mathcal{N}$ denotes the region in \mbox{$s$-$t$} parameter space that maps to the region on the manifold lying inside the search sphere of radius $r_s$ centered at the point $\bar{v}$.
The area element $dS$ is given by the formula $dS=\norm{\tfrac{\partial x}{\partial s}(s,t) \times \tfrac{\partial x}{\partial t}(s,t)}_2 ds \, dt$, where $s$ and $t$ are the manifold parameters, and $\norm{\cdot}_2$ denotes the standard norm.
The integral in the denominator is the surface area of the manifold inside the search sphere, which serves as a normalization constant.

To elucidate the relationship between the continuous version of $P(\bar{v}, \mathcal{N})$ and the curvature of the underlying manifold, we will focus on the case where the underlying manifold is a quadratic surface.
Given a quadratic surface in $\mathbb{R}^3$ and a point on the surface, there is always a translation of the surface and an orthonormal change of variables corresponding to a rotation so that we may assume without a loss of generality that the point under consideration is the origin, and the normal to the surface is aligned with the $x_3$ axis.
Such a quadratic surface has the form
\begin{equation}
\label{eqn:qdrtc_srfc}
0 = 
x^T M x + x_3
\end{equation}
where
\begin{equation*}
M \equiv
\begin{bmatrix}
M_{11} & 0 & \frac{1}{2} M_{13} \\
0 & M_{22} & \frac{1}{2} M_{23} \\
\frac{1}{2} M_{13} & \frac{1}{2} M_{23} & M_{33}
\end{bmatrix}
\text{ and }
x =
\begin{bmatrix} x_1 \\ x_2 \\ x_3 \end{bmatrix}
\end{equation*}
An expansion in the search sphere radius $r_s$ of the local direction covariance matrix (\ref{eqn:exct_emprcl_drctn_mtrx}) of the quadratic surface (\ref{eqn:qdrtc_srfc}) at the origin is
\begin{equation}
\label{eqn:qdrtc_lcl_drctn_cvrnc}
P(0, \mathcal{N})
=
\begin{bmatrix}
\frac{1}{2} + a_{11} r_s^2 & 0 & a_{13} r_s^2 \\
0 & \frac{1}{2} + a_{22} r_s^2 & a_{23} r_s^2 \\
a_{13}r_s^2 & a_{23} r_s^2 & a_{33} r_s^2
\end{bmatrix}
+ \cdots
\end{equation}
where
\begin{gather*}
\begin{aligned}
a_{11} = -\frac{1}{32}\bigpars{(M_{11}+M_{22})^2 + 4 M_{11}^2}
& &
a_{22} = -\frac{1}{32}\bigpars{(M_{11}+M_{22})^2 + 4 M_{22}^2}
\end{aligned}
\\
\begin{aligned}
a_{13} = \frac{1}{16} M_{13}(3 M_{11} + M_{22})
& &
a_{23} = \frac{1}{16} M_{23}(M_{11} + 3 M_{22})
\end{aligned}
\\
a_{33} = \frac{1}{16}( 3 M_{11}^2 + 2 M_{11} M_{22} + 3 M_{22}^2)
\end{gather*}
For the remainder of the discussion, we refer to the local direction covariance matrix to leading order in (\ref{eqn:qdrtc_lcl_drctn_cvrnc}) as $P$.
Note that $P$ is diagonal when $r_s=0$, and its eigenvalues are $1/2$ and $0$ with algebraic multiplicities two and one.
When $r_s \ne 0$, the eigenvalues of $P$ are perturbations off of $1/2$ and $0$ that depend on the curvature of the quadratic surface and the size of the search radius $r_s$.
If a dominant eigenvalue of $P$ has an expansion of the form $\lambda = 1/2 + c_1 r_s + c_2 r_s^2 + \cdots$, then to leading order $\lambda = 1/2 + a_{11}r_s^2$ and $\lambda = 1/2 + a_{22} r_s^2$.
By the same argument, the weak eigenvalue of $P$ to leading order is $\lambda = a_{33} r_s^2$.
Both these approximations to the dominant eigenvalues satisfy the characteristic polynomial of $P$ with a residual of $\mathcal{O}(r_s^4)$, where the characteristic polynomial is computed from the leading order approximation \eqref{eqn:qdrtc_lcl_drctn_cvrnc}.
The signed curvature at the origin\footnote{The signed curvature at $x$ of the plane curve $(x, y(x))$ is $\kappa=y''(x)/(1+y'(x)^2)^{3/2}$} of the curve formed by the intersection of the quadratic surface \eqref{eqn:qdrtc_srfc} and the plane perpendicular to the tangent space containing the vector $[ \cos(\theta) \phantom{+} \sin(\theta) \phantom{+} 0]^T$ is $\kappa(\theta) = -2(M_{11} \cos(\theta)^2 + M_{22} \sin(\theta)^2)$.
It follows that the maximum and minimum curvatures are $\kappa_1 \equiv -2 M_{11}$ and $\kappa_2 \equiv -2 M_{22}$.
The two directions in the tangent space associated with the curvatures $\kappa_1$ and $\kappa_2$ are aligned with the $x_1$ and $x_2$ axes.
Let $\bar{\kappa}=(\kappa_1 + \kappa_2)/2$ denote the mean curvature at the origin.
The eigenvalues of $P$ in terms of the mean, maximum, minimum curvatures, and search radius are
\begin{equation}
\begin{aligned}
\lambda_1(P) &= \frac{1}{2} - \frac{1}{32} (\bar{\kappa}^2 + \kappa_1^2)r_s^2 + \cdots
\\
\lambda_2(P) &= \frac{1}{2} - \frac{1}{32} (\bar{\kappa}^2 + \kappa_2^2)r_s^2 + \cdots
\\
\lambda_3(P) &= \frac{1}{32} ( 2\bar{\kappa}^2 + \kappa_1^2 + \kappa_2^2 ) r_s^2 + \cdots
\end{aligned}
\label{eqn:qdrtc_lcl_cvrnc_evals}
\end{equation}
When the search radius $r_s$ is zero, the eigenvectors associated with dominant eigendirections of $P$ are the first two columns of the identity matrix, and they form a basis for the tangent space through the origin of the quadratic surface.
When $r_s \ne 0$, then in general the eigenvectors of $P$ no longer span the tangent space, but they span a plane that nearly coincides with the tangent space.
To calculate an approximate eigenvector $u^1$ of $P$ associated with $\lambda_1(P)$ \eqref{eqn:qdrtc_lcl_cvrnc_evals}, we may safely assume that $u^1$ has been scaled so that the first entry is $u^1_1 = 1$, and the other entries have Taylor expansions in the search radius $r_s$.
The coefficients in the Taylor expansions can be calculated by the standard technique of zeroing the coefficients of the leading order terms in $r_s$ of $\norm{(P - \lambda_1(P) I) u^1}_2^2$.

The approximate eigenvectors of $P$ are
\begin{equation}
\begin{aligned}
u^1 &= e^1 + 2 a_{13}r_s^2 e^3 + \cdots \\
u^2 &= e^2 + 2 a_{23}r_s^2 e^3 + \cdots 
\label{eqn:qdrtc_P_dmnt_apprx_evecs}
\end{aligned}
\end{equation}
where $e^i$ denotes the $i^{\text{th}}$ column of the identity matrix.
These formulas indicate that $u^1$ and $u^2$ are nearly aligned with the directions in which the curvature of the quadratic surface is maximized and minimized.
For both approximate eigenvectors, $\norm{(P - \lambda_i(P)I)u^i}_2 = \mathcal{O}(r_s^4)$ after $u^i$ has been normalized in the two-norm.

%%%%%%%%%%%%%%%%%%%%%%%%%%%%%%%%%%%%%%%%%%%%%%%%%%%%%%%%%%%%%%%%%%%%%%%%%%%%%%%%%%%%%%%%%%%%%%%%%%%%%%%%%%
\subsection{Motivation for the vertex placement algorithm}
\label{sec:cntns_vrtx_plcmnt}

We can adapt the constrained minimization problem \eqref{eqn:cnsrtnd_optmzn_prblm} that underlies the new vertex placement algorithm of \S \ref{sec:R3_gnrtng_new_trngl} to a continuous setting by treating the underlying manifold as the data.
In this setting, the vertices $v^1$ and $v^2$ are points on the manifold, and the matrices $Q^1$ and $Q^2$ are computed from the continuous version of the local direction covariance matrix \eqref{eqn:exct_emprcl_drctn_mtrx}.
The main requirement of the new vertex placement algorithm is that the vertex must not be too far from the underlying surface.
Otherwise, the resulting triangle cannot possibly fit the data in a meaningful way.
Here, too far means a large distance as measured by either of the two metrics induced by $Q^1$ and $Q^2$. 
These metrics are built from curvature information local to the active edge vertices, which can be estimated in the discrete setting where we just have points on the underlying surface, but not the surface itself.
With this in mind, the initial vertex placement algorithm may be interpreted as simultaneously minimizing the distances from the new vertex to the underlying surface as measured by the induced metrics associated with each active edge vertex.
There appears to be no reason to favor one metric over the over, and the second constraint in \eqref{eqn:opmtzn_cnstrnts} ensures both are given equal weight in the simultaneous minimization.

If one of the front edge vertices is the origin, then the objective function \eqref{eqn:opmtzn_objctv}
is $x^T Q x$, where $Q = (P + \mu I)^{-1}$ and $\mu$ is a user supplied small parameter.
In the context of the minimization problem that underlies the initial vertex placement algorithm, it is helpful to think of $x^T Q x$ as a penalty function that penalizes points based on their distance to the underlying surface.
Under this interpretation, a point $x$ that is not near the underlying surface is penalized in the sense that $x^T Q x$ is large.
For a general point in space, the size of the penalty is determined by the eigen structure of $Q$, which is intimately related to the curvature and tangent space of the underlying surface.
In the case of the quadratic surface \eqref{eqn:qdrtc_srfc}, the exact tangent space through the origin is the $x_1$-$x_2$ plane, and we define the \emph{approximate tangent space} to be the subspace spanned by  $u^1$ and $u^2$ in \eqref{eqn:qdrtc_P_dmnt_apprx_evecs}.
$P$ and $Q$ share the same eigenvectors, so two of $Q$'s eigenvectors span the approximate tangent space, and its third eigenvector is orthogonal to the the approximate tangent space.
The largest eigenvalue of $Q$ is $\lambda_{\text{max}}(Q) = (\mu + \lambda_{\text{min}}(P))^{-1} = \mathcal{O}((\mu + r_s^2)^{-1})$, which is much larger than the other two eigenvalues assuming the user has not chosen an overly large value of $\mu$.
Consequently, if $x$ has a non-negligible component in the direction perpendicular to the approximate tangent space, then the objective function severely penalizes $x$ in the sense that the value $x^T Q x$ is dominated by $\lambda_{\text{max}}(Q)$.
The objective function should exhibit this behavior since points that are not near the approximate tangent space cannot be near the underlying surface. 
The remaining two eigendirections of $Q$ lie in the approximate tangent space, and these eigendirections are nearly aligned with the directions of maximum and minimum curvature of the quadratic surface as described in \S \ref{sec:lcl_cvrnc_mtrx_intrptn}.
The eigenvalues of $Q$ in these directions are
\begin{equation}
\begin{aligned}
\lambda_1(Q) &= \Bigpars{\frac{1}{2} + \mu -\frac{1}{32}(\bar{\kappa}^2 + \kappa_1^2) r_s^2 + \cdots }^{-1}
\\
\lambda_2(Q) &= \Bigpars{\frac{1}{2} + \mu -\frac{1}{32}(\bar{\kappa}^2 + \kappa_2^2) r_s^2 + \cdots }^{-1}
\end{aligned}
\label{eqn:qdrtc_Q_weak_evals}
\end{equation}
There is enough freedom in the derivation of the quadratic surface equation \eqref{eqn:qdrtc_srfc} that we may assume for the sake of concreteness that $\abs{\kappa_1} \le \abs{\kappa_2}$, which then imposes $\lambda_1(Q) \le \lambda_2(Q)$.
Thus, for two points with equal Euclidean norm that both lie in the approximate tangent space, the objective function penalizes the point with the greater component in the $\lambda_2(Q)$ direction more than it penalizes the other point.
The eigendirection associated with $\lambda_2(Q)$ nearly coincides with the direction of greatest curvature in magnitude, which is $\abs{\kappa_2}$.
It is in this direction that the quadratic surface curves away from the exact tangent space most rapidly.
The approximate and exact tangent spaces are nearly aligned, so points with a larger component in the direction of the maximum magnitude curvature will be farther from the underlying quadratic surface than points with the same Euclidean norm, but a smaller component in the direction of maximum magnitude curvature.
The curvature of the surface is encoded in $Q$ through its eigenvalues \eqref{eqn:qdrtc_Q_weak_evals}, which allows the objective function to penalize points in the approximate tangent space that are far from the underlying surface by using the local curvature information as a proxy.

%%%%%%%%%%%%%%%%%%%%%%%%%%%%%%%%%%%%%%%%%%%%%%%%%%%%%%%%%%%%%%%%%%%%%%%%%%%%%%%%%%%%%%%%%%%%%%%%%%%%%%%
%% START: motivating  the constrained minimization problem with an ellipsoid. Subsumed by more general explanation
%%%%%%%%%%%%%%%%%%%%%%%%%%%%%%%%%%%%%%%%%%%%%%%%%%%%%%%%%%%%%%%%%%%%%%%%%%%%%%%%%%%%%%%%%%%%%%%%%%%%%%%

%%%%%%%%%%%%%%%%%%%%%%%%%%%%%%%%%%%%%%%%%%%%%%%%%%%%%%%%%%%%%%%%%%%%%%%%%%%%%
\subsection{Conditions under which the minimization problem has a solution}
Although the constrained optimization problem (\ref{eqn:cnsrtnd_optmzn_prblm}) is not guaranteed to have a solution, we can derive mild conditions under which it does so that the algorithm user can be confident that the advancing front stage of the algorithm is robust.
We focus on the case where  we are triangulating data in $\mathbb{R}^N$.
We omit the third constraint in (\ref{eqn:opmtzn_cnstrnts}) whose purpose is to make the minimizer unique, so that the simplified constraint set is the intersection of the constraint sphere and the surface that describes the isosceles constraint.
Without a loss of generality, we may assume that $v^2 = -v^1$, so that the midpoint $v^m$ of the two active edge vertices falls on the origin.
First, consider the case where $Q^1 =  Q^2$, then the isosceles constraint reduces to $0 = v^T Q^1 v^1$ which describes a plane through the origin.
The constraint set in this case is the intersection of the plane with a sphere centered at the origin.
Therefore, the constraint set is a non-empty compact subset of $\mathbb{R}^N$, which guarantees that a minimizer of the constrained optimization problem exists.

Now, consider the more realistic case where $Q^1 \ne Q^2$, but the tangent spaces through $v^1$ and $v^2$ are nearly aligned in the sense that their canonical angles\footnote{Let $U$ and $V$ denote orthonormal matrices whose columns span two subspaces of equal dimension.
The cosines of the canonical subspace angles are the singular values of $U^TV$  \cite[p. 73]{0898714141}.}
are nearly zero, and the eigenvalues of $Q^1$ and $Q^2$ are nearly equal after they have been sorted by magnitude.
The isosceles constraint reduces to
\begin{equation}
\label{eqn:iscles_cnstrnt_prtrbtn}
0 = 4 v^T Q^1 v^1 + 2 v^T E_Q v^1 + v^T E_Q v
\end{equation}
where $E_Q \equiv Q^2 - Q^1$.
The first term in equation (\ref{eqn:iscles_cnstrnt_prtrbtn}) describes the plane from the isosceles constraint in the $Q^1=Q^2$ case, and the other terms can be interpreted as perturbations that introduce curvature into the isosceles constraint surface.
By continuity, if $E_Q$ is sufficiently small in norm, then the surface that satisfies the isosceles constraint nearly coincides with the plane that describes the isosceles constraint in the $Q^1=Q^2$ case.
It follows that if $E_Q$ is sufficiently small, then the constraint set described by the intersection of the constraint sphere with the isosceles constraint surface is a compact subset of $\mathbb{R}^N$, and therefore a minimizer to the constrained minimization problem exists.
We now show that the magnitude of $E_Q$ is controlled by the alignment of the two tangent spaces, and the difference in the eigenvalues of $Q^1$ and $Q^2$.
If the data points reside on a smooth manifold, $v^1$ and $v^2$ lie near each other on the manifold, and if there are a sufficient number of data points in the neighborhoods of $v^1$ and $v^2$, then the alignment of the two tangent spaces and difference in eigenvalues of $Q^1$ and $Q^2$ are governed by the curvature of the underlying manifold.
Furthermore, as $v^1$ and $v^2$ approach each other and the density of data in their respective neighborhoods increases, the tangent spaces approach perfect alignment and the differences between eigenvalues approach zero. 
The symmetric matrices $Q^1$ and $Q^2$ can be diagonalized by orthonormal transformations, so let $Q^1 = U^1 \Lambda^1 (U^1)^T$ and $Q^2 = U^2 \Lambda^2 (U^2)^T$ where $U^i$ is an orthonormal matrix and $\Lambda^i$ is diagonal, and define $E_U \equiv (U^1)^T U^2 - I$ and $E_{\Lambda} \equiv \Lambda_2 - \Lambda_1$.
%Since the tangent spaces associated with $v^1$ and $v^2$ are nearly aligned, we may assume without a loss 
We may assume without a loss of generality that $\norm{E_U}_2$ is on the order of $\max_i(1-\cos(\theta_i))$, where $\theta_i$ denotes one of the canonical angles between the two tangent spaces, or the canonical angle between the orthogonal complements of the tangent spaces.
Thus, $\norm{E_U}$ approaches zero as the two tangent spaces approach perfect alignment.
By a standard argument,
\begin{equation*}
\norm{E_Q}_2 \le \norm{\Lambda^1 - \Lambda^2}_2 + 2 \norm{\Lambda^2}_2 \norm{E_u}_2 + \norm{\Lambda^2}_2 \norm{E_u}_2^2
\end{equation*}
which implies that $\norm{E_Q}_2$ is bounded by $\norm{ \Lambda_2 - \Lambda_1}_2$ and $\norm{E_U}_2$.

%%%%%%%%%%%%%%%%%%%%%%%%%%%%%%%%%%%%%%%%%%%%%%%%%%%%%%%%%%%%%%%%%%%%%%%%%%%%%
\section{Results}
\label{sec:results}
We tested the SNPCA algorithm by running it on data sets that fall near the surface of a sphere, torus, swiss roll, and creased sheet.
All these manifolds can be parameterized by two variables, and we embedded the data in a fifty dimensional space as follows.
To generate points in $\mathbb{R}^3$ that are distributed uniformly over the manifold with respect to surface area, we sampled parameter values from the distribution whose probability density function is the manifold area element normalized by the manifold surface area.
For example, the sphere of radius $a$ has surface area $4 \pi a^2$ and is parameterized by $x(u,v) = a (\cos(u) \sin(v), \sin(u) \sin(v), \cos(v))$ where $0 \le u \le 2 \pi$ and $0 \le v \le \pi$.
The area element is $dS = a^2 \sin(v) du \, dv$, so the distribution from which we sample values of $u$ and $v$ for the sphere is $f(u,v) = \sin(v)/4 \pi$.
We then generated points in $\mathbb{R}^3$ on the manifold from these random parameter values.
Finally, we randomly generated three orthogonal unit vectors in $\mathbb{R}^{50}$, and then for each point on the manifold, we generated a data point in $\mathbb{R}^{50}$ by forming the linear combination of three three orthonormal vectors with the manifold point coordinates as weights.
The vertex coordinates of the final triangulation are in $\mathbb{R}^{50}$, so to produce the illustrations of the final triangulation, we undid the orthogonal transformation so that the vertex coordinates were in $\mathbb{R}^3$.

We chose these data sets based on attributes of the manifold that underlies the data, namely the curvature and smoothness of the manifold, and the presence of a boundary.
The creased sheet is smooth with no curvature in every neighborhood away from the crease location, and has a boundary.
The sphere is smooth, has uniform curvature, and no boundary
The swiss roll is smooth, has varying curvature, and a boundary.
The torus is smooth, has varying curvature, and no boundary.

For each data set, we have included illustrations of the points in the data set, the output of the advancing front stage, and the final triangulation produced by the seam sewing stage.
The error tables list the maximum error, the average error, and the RMS error associated with each triangulation generated by SNPCA.
The error associated with the $i^{\text{th}}$ data point, denoted by $d_i$, is the shortest distance from the data point to a point on the triangulation.

%%%%%%%%%%%%%%%%%%%%%%%%%%%%%%%%%%%%%%%%%%%%%%%%%%%%%%%%%%%%%%%%%%%%%%%%%%%%%%%%%%%%%%%%%%%%%%%%%%%%%%%%%
%\clearpage
%\clearpage
%\subsection{Sphere}
%The surface underlying this data set is a sphere of radius $4$ in $\mathbb{R}^3$.

\begin{figure}[h!]
\begin{center}
\complextable{$1372$}{$689$}{$0.5$}\\
\errtable{$0.082777$}{$0.013134$}{$0.014518$}{$10000$}
\subfloat[Data]{\includegraphics[scale=.45]{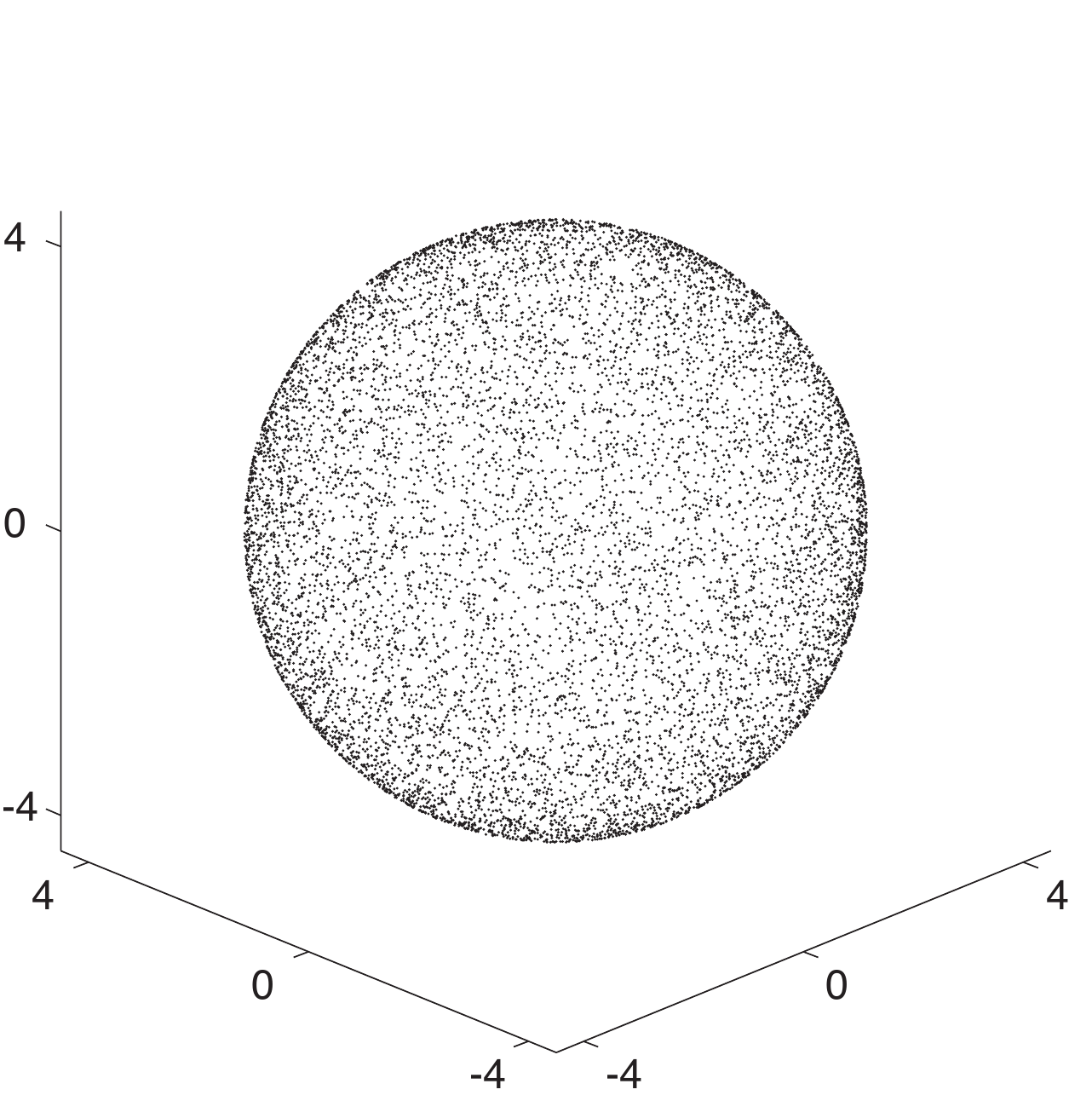}}
\subfloat[SNPCA seams]{\includegraphics[scale=.45]{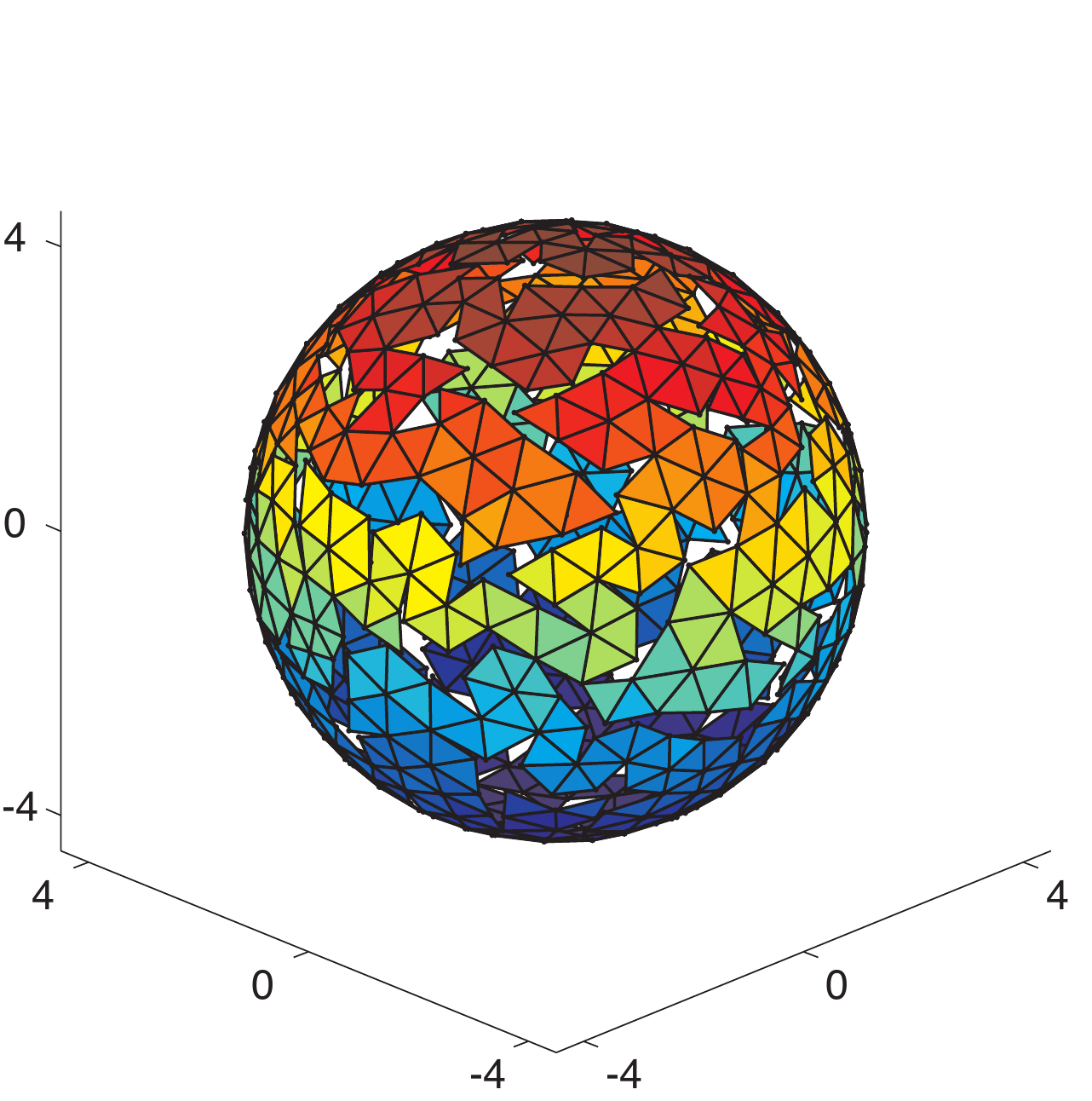}}
\\
\subfloat[SNPCA sewn]{\includegraphics[scale=.45]{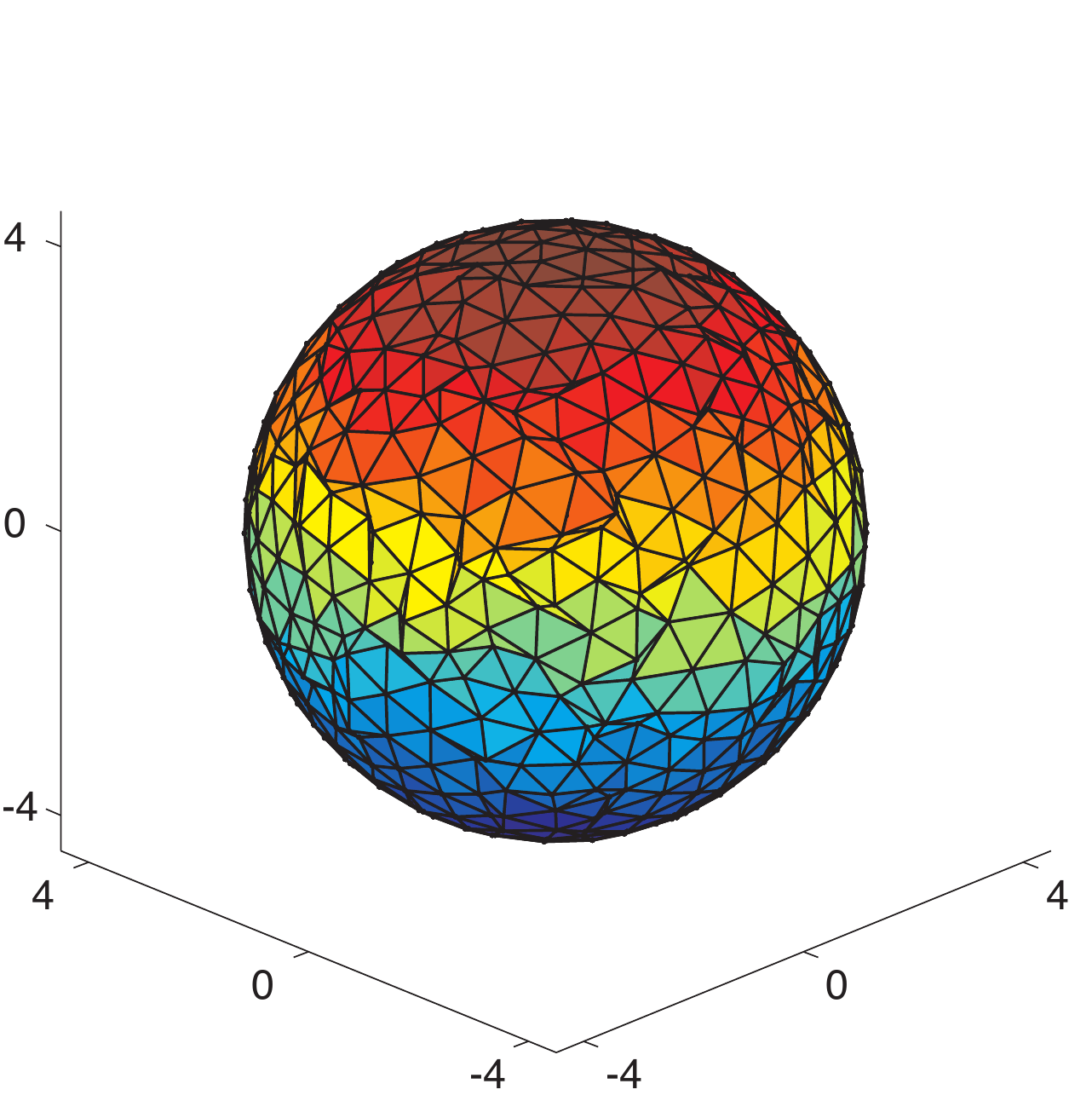}}
%\subfloat[PCA]{}%\includegraphics[scale=.45]{sphere_PCA_gray.pdf}}
\label{fig:sphere}
%\caption{{\bf Sphere}}
\caption{A $2$-sphere of radius $4$ embedded in $\mathbb{R}^{50}$.}
\end{center}
\end{figure}

%\clearpage

%%%%%%%%%%%%%%%%%%%%%%%%%%%%%%%%%%%%%%%%%%%%%%%%%%%%%%%%%%%%%%%%%%%%%%%%%%%%%%%%%%%%%%%%%%%%%%%%%%%%%%%%%
%\subsection{Torus}
%The surface underlying this data set is a torus with radii of $4$ and $1$.
\begin{figure}[h!]
\begin{center}
\complextable{$681$}{$342$}{$0.75$}\\
\errtable{$0.10804$}{$0.034369$}{$0.040669$}{$10000$}
\subfloat[Data]{\includegraphics[scale=.35]{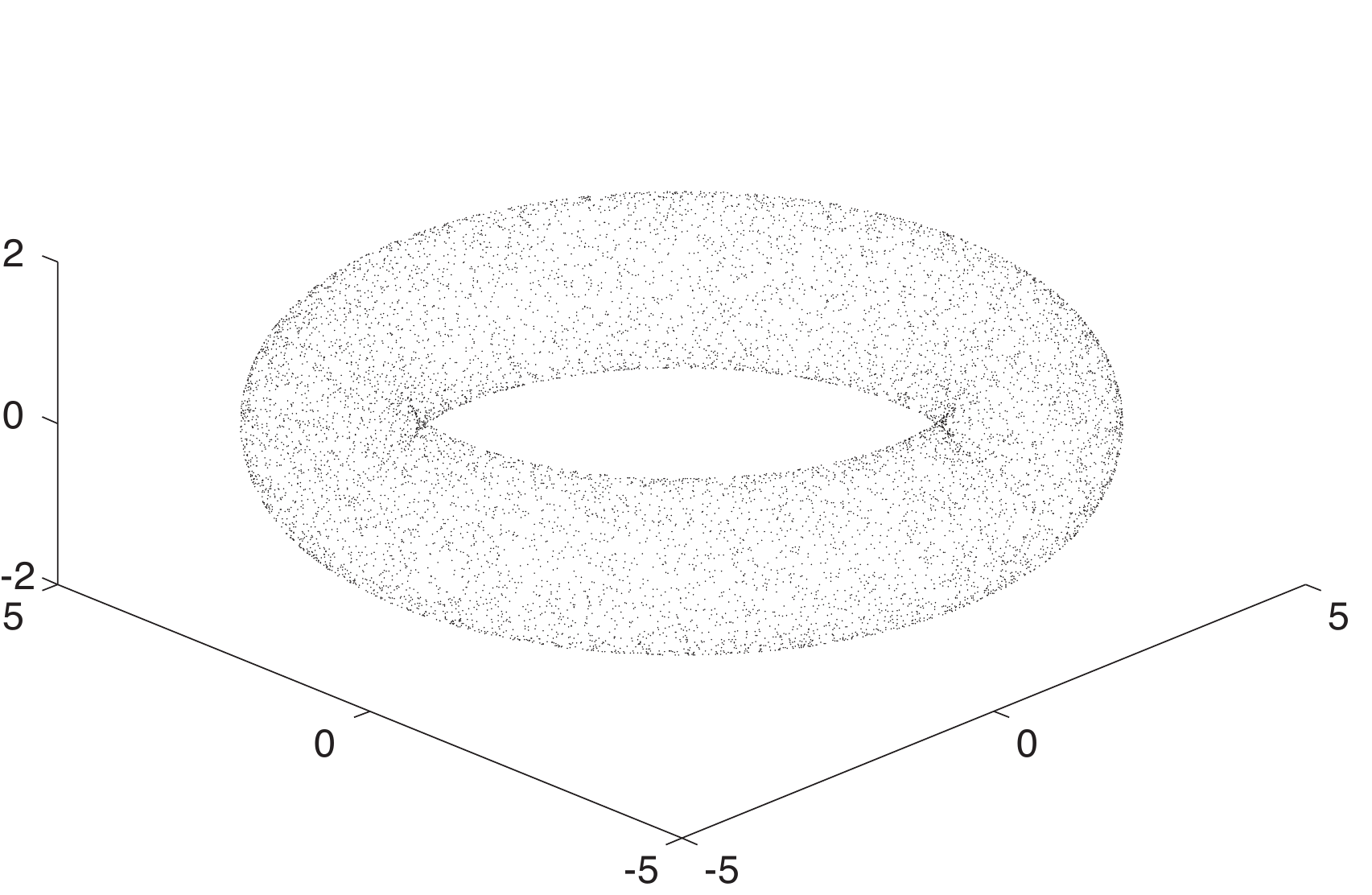}}
\subfloat[SNPCA seams]{\includegraphics[scale=.35]{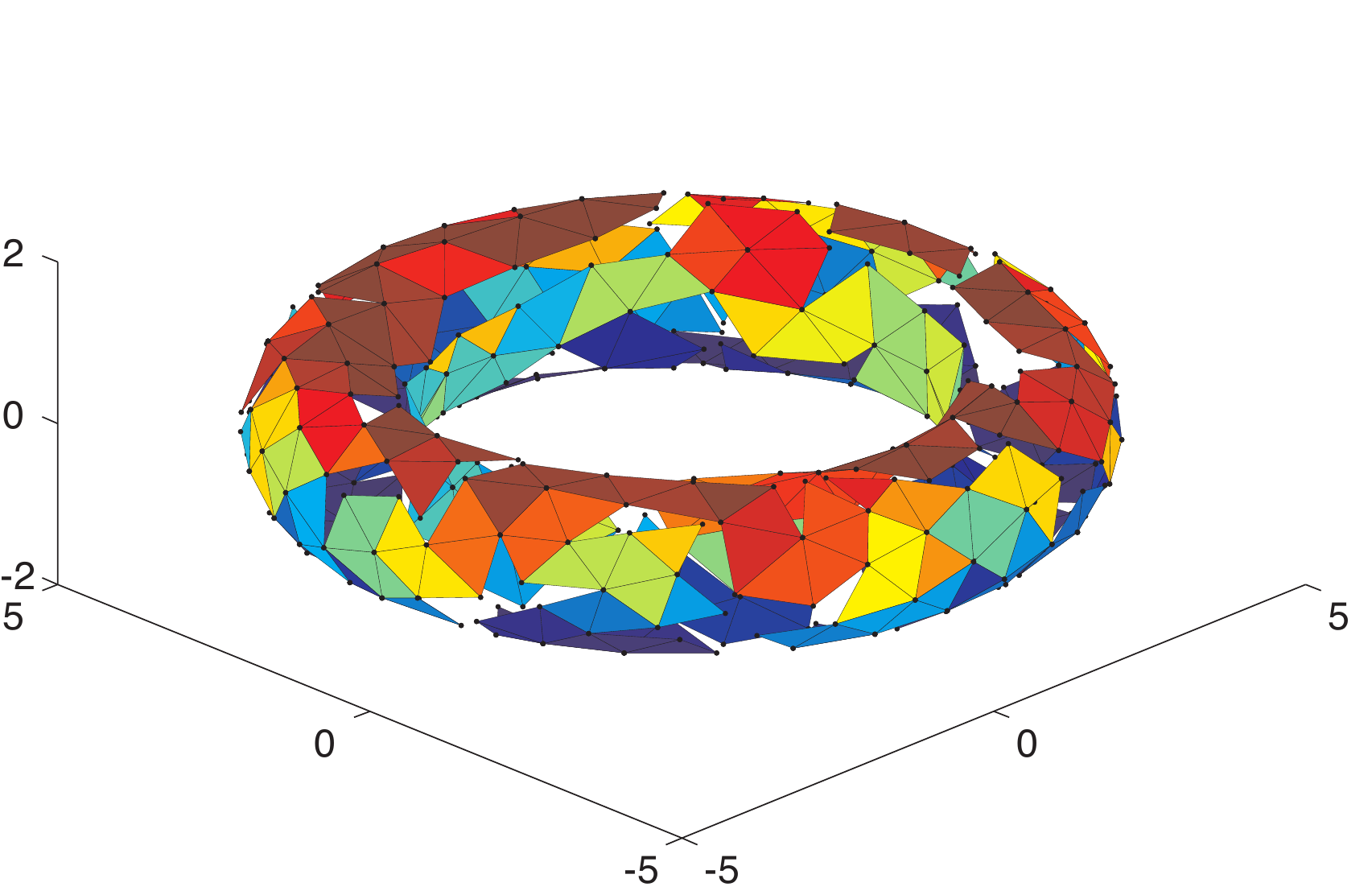}}
\\
\subfloat[SNPCA sewn]{\includegraphics[scale=.35]{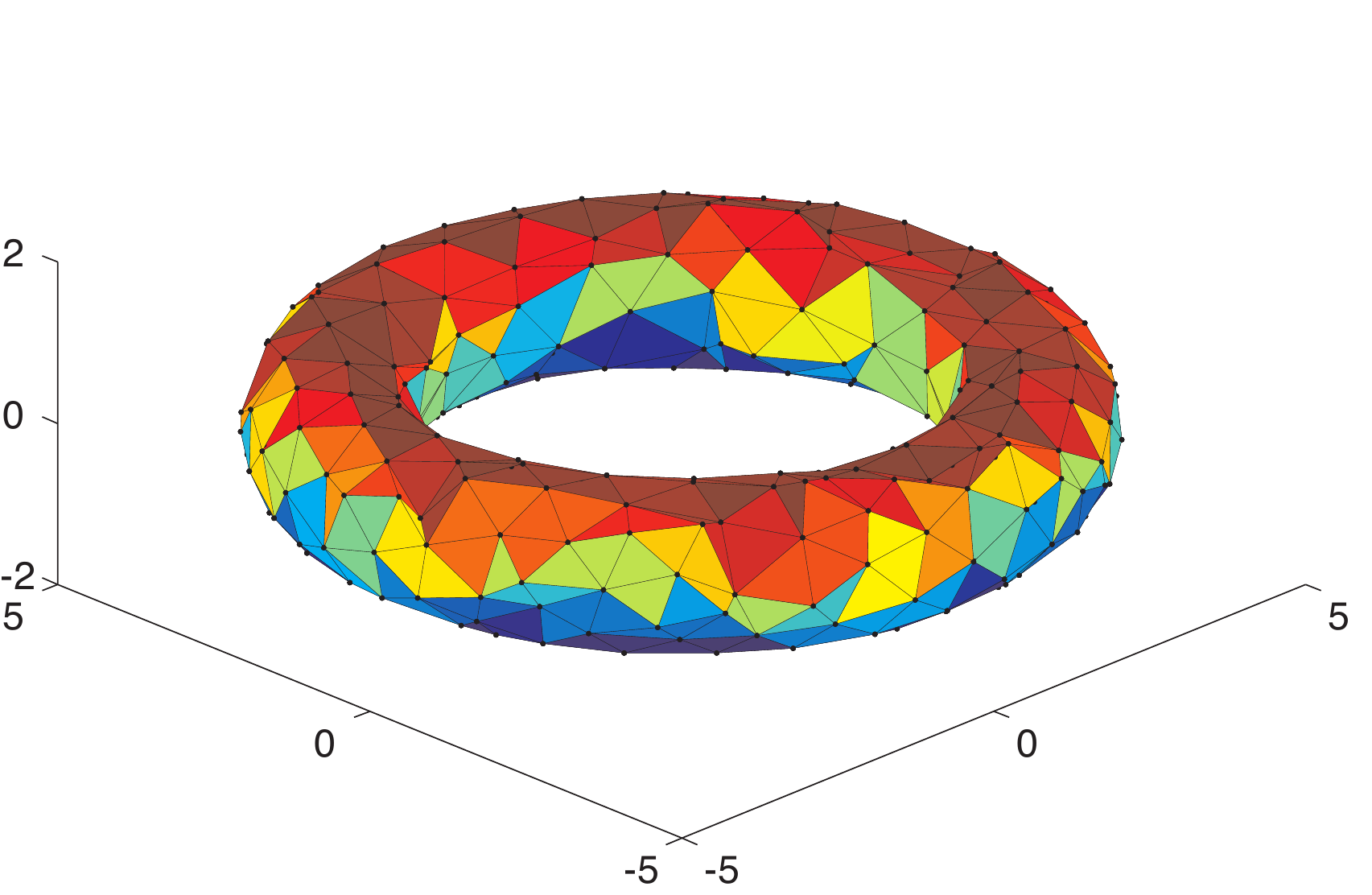}}
%\subfloat[PCA]{}%\includegraphics[scale=.35]{torus_short_PCA_gray.pdf}}
\label{fig:torus}
\caption{A $2$-torus with radii $4$ and $1$ embedded in $\mathbb{R}^{50}$.}
\end{center}
\end{figure}

%\clearpage

%%%%%%%%%%%%%%%%%%%%%%%%%%%%%%%%%%%%%%%%%%%%%%%%%%%%%%%%%%%%%%%%%%%%%%%%%%%%%%%%%%%%%%%%%%%%%%%%%%%%%%%%%
%\subsection{Torus with noise}

\begin{figure}[h!]
\begin{center}
\complextable{$1245$}{$639$}{$0.5$}\\
\errtable{$0.65358$}{$0.088305$}{$0.11857$}{$10000$}
\subfloat[Data]{\includegraphics[scale=.35]{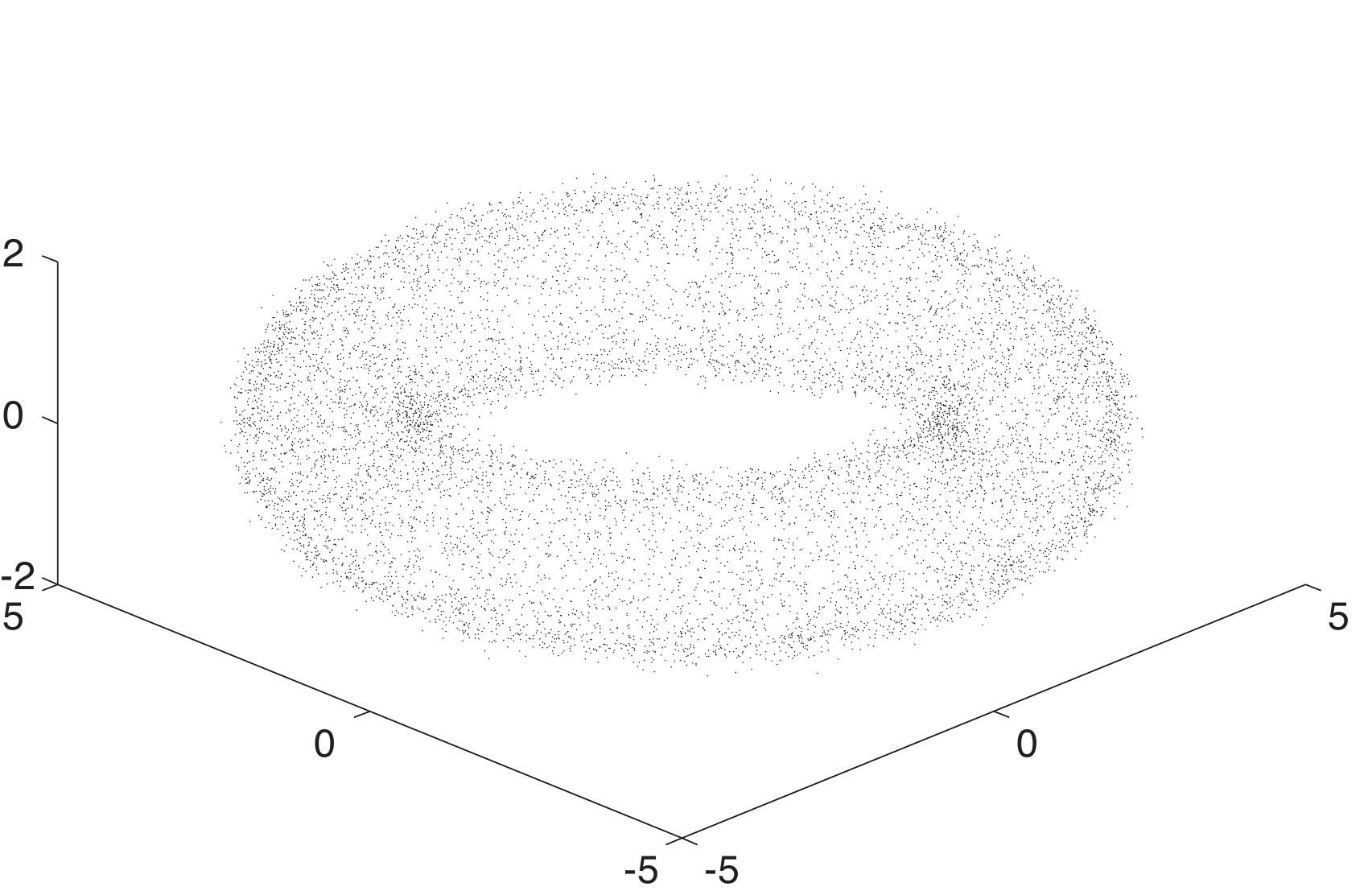}}
\subfloat[SNPCA seams]{\includegraphics[scale=.35]{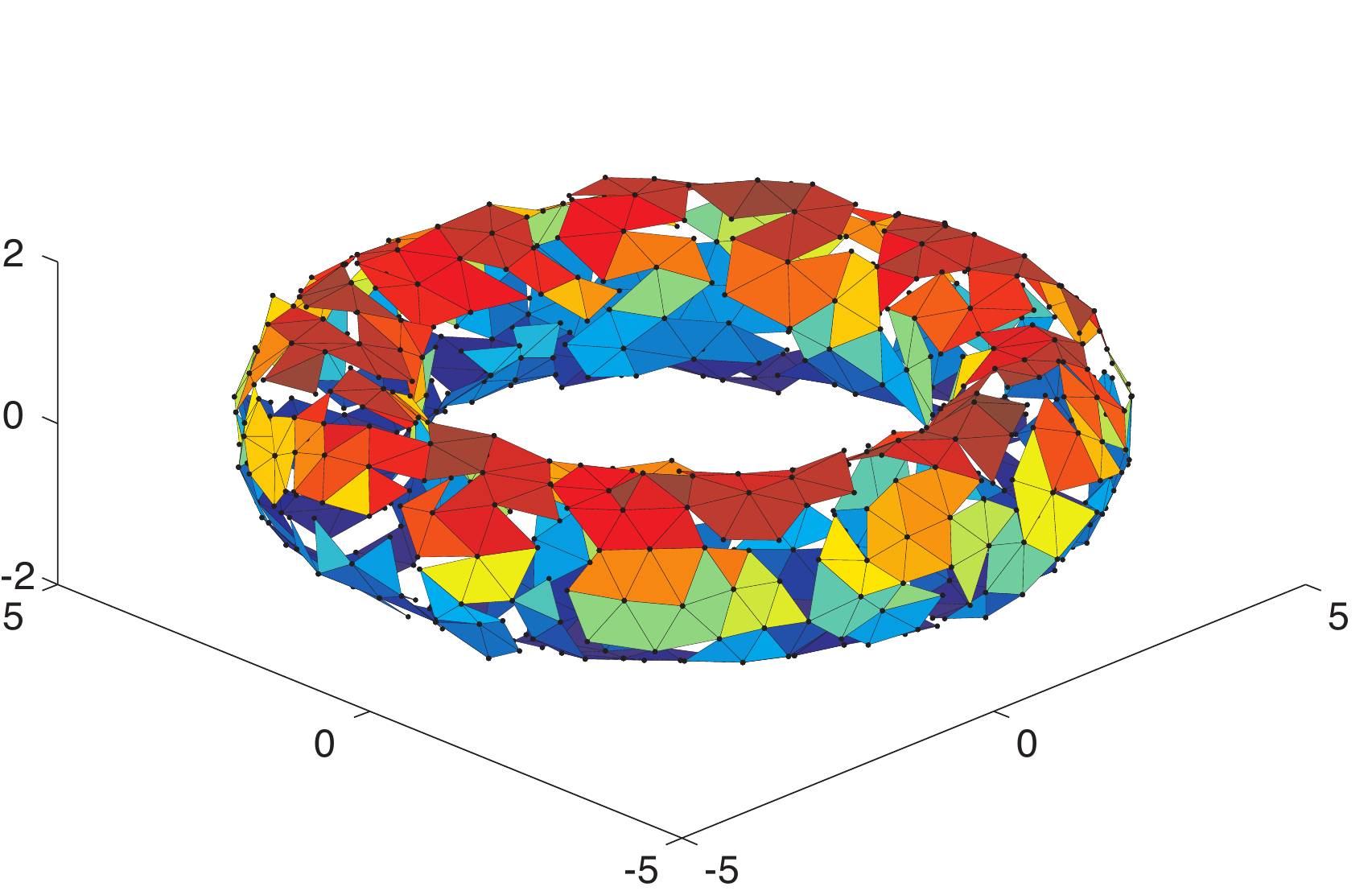}}
\\
\subfloat[SNPCA sewn]{\includegraphics[scale=.35]{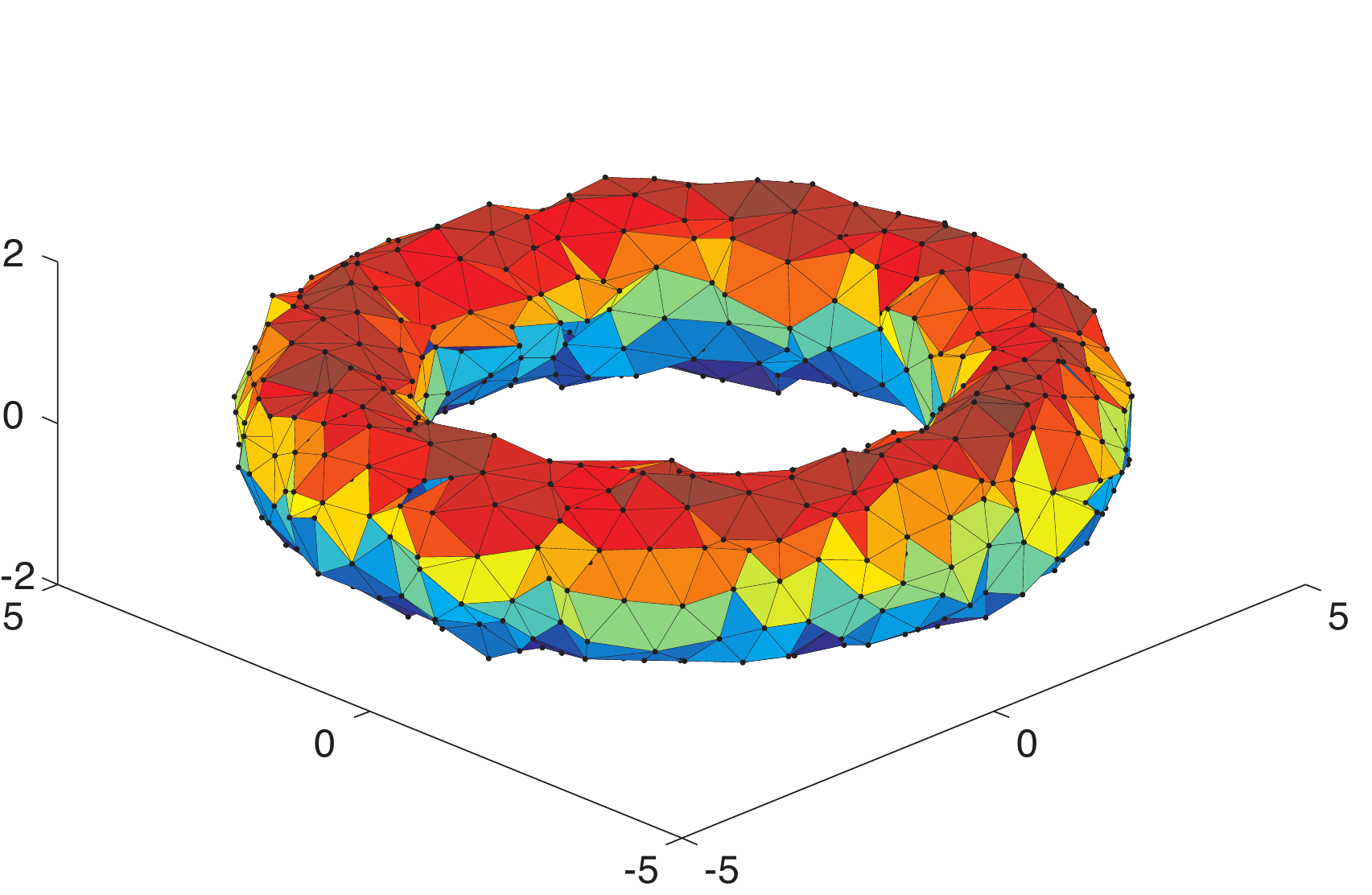}}
%\subfloat[PCA]{}%\includegraphics[scale=.35]{torus_noise_PCA_gray.pdf}}
\label{fig:torus_noise}
\caption{The surface underlying this data set is a torus with radii of $4$ and $1$, and we added error to the surface data points before generating the set of data vectors.
We generated the error vectors in $\mathbb{R}^3$ by drawing each coordinate from the standard normal distribution, and then scaling the error vector by a factor of $.01$.
We then generated the set of data vectors by our standard orthogonal transformation.
The final triangulation is not \emph{watertight}, meaning it has front edges.
}
\end{center}
\end{figure}

%\clearpage

%%%%%%%%%%%%%%%%%%%%%%%%%%%%%%%%%%%%%%%%%%%%%%%%%%%%%%%%%%%%%%%%%%%%%%%%%%%%%%%%%%%%%%%%%%%%%%%%%%%%%%%%%
%\subsection{Swiss roll}

\begin{figure}[h!]
\begin{center}
\complextable{$698$}{$399$}{$0.75$}\\
\errtable{$0.59827$}{$0.015952$}{$0.031548$}{$10000$}
\subfloat[Data]{\includegraphics[scale=.40]{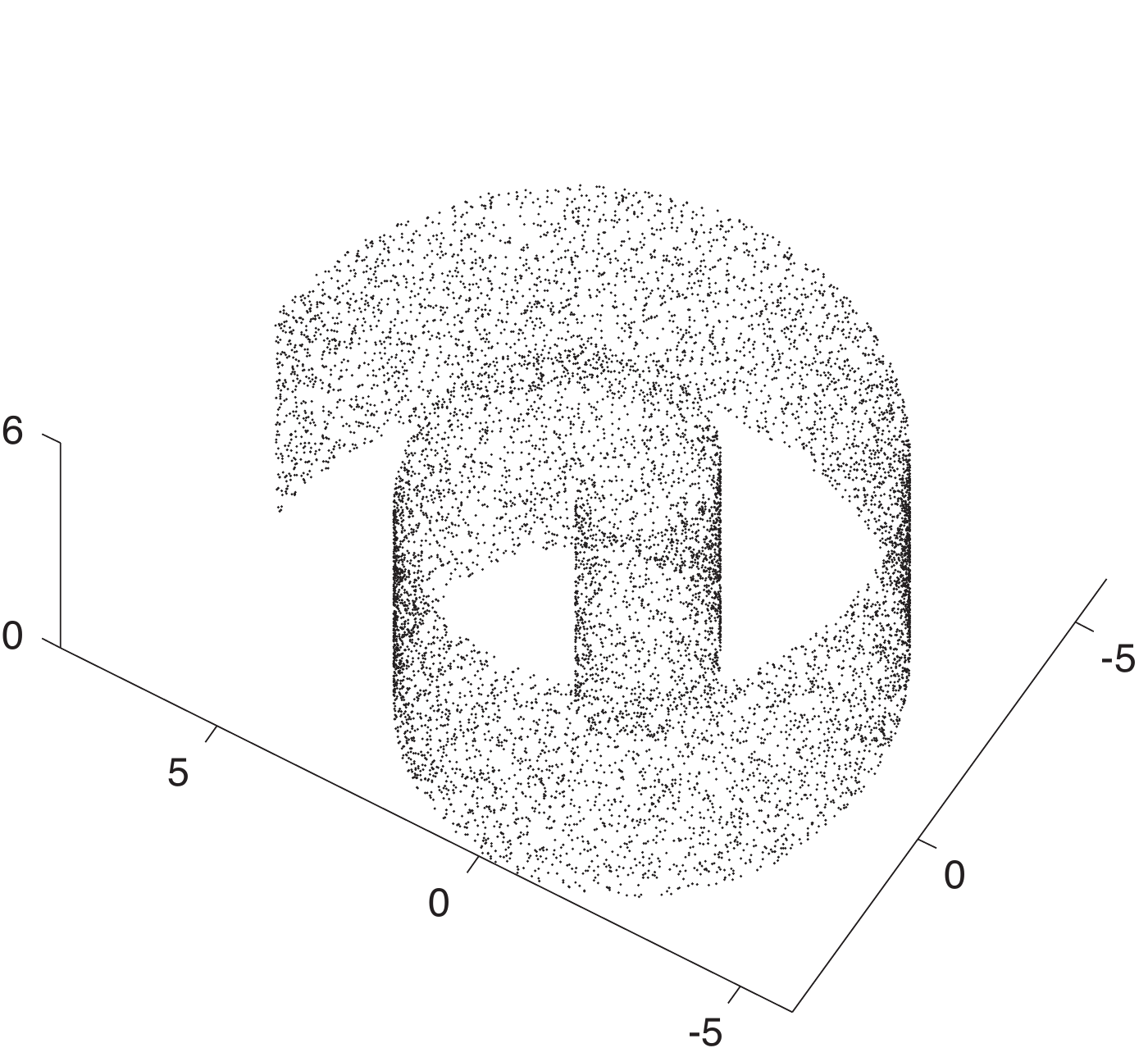}}
\subfloat[SNPCA seams]{\includegraphics[scale=.40]{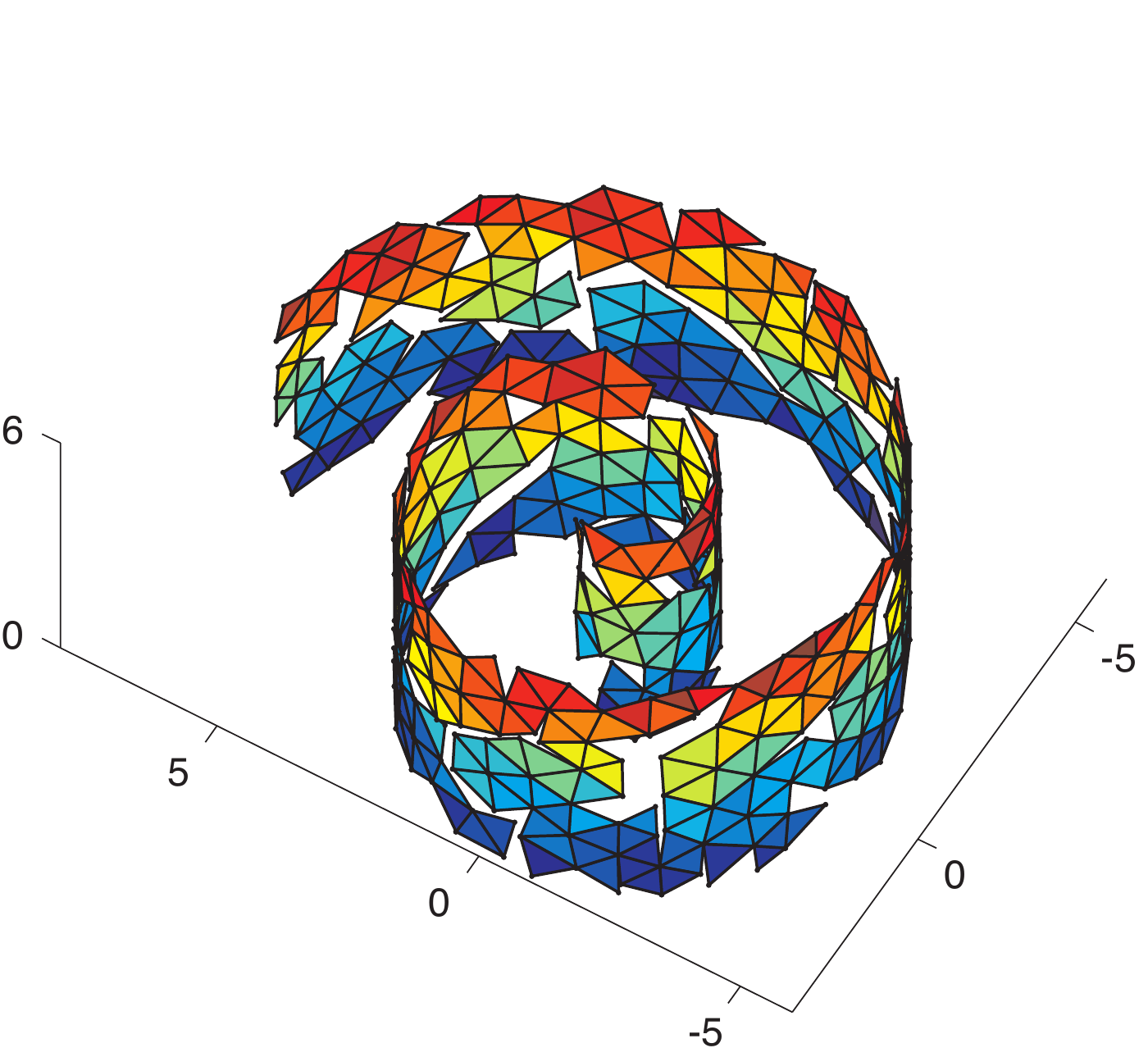}}
\\
\subfloat[SNPCA sewn]{\includegraphics[scale=.40]{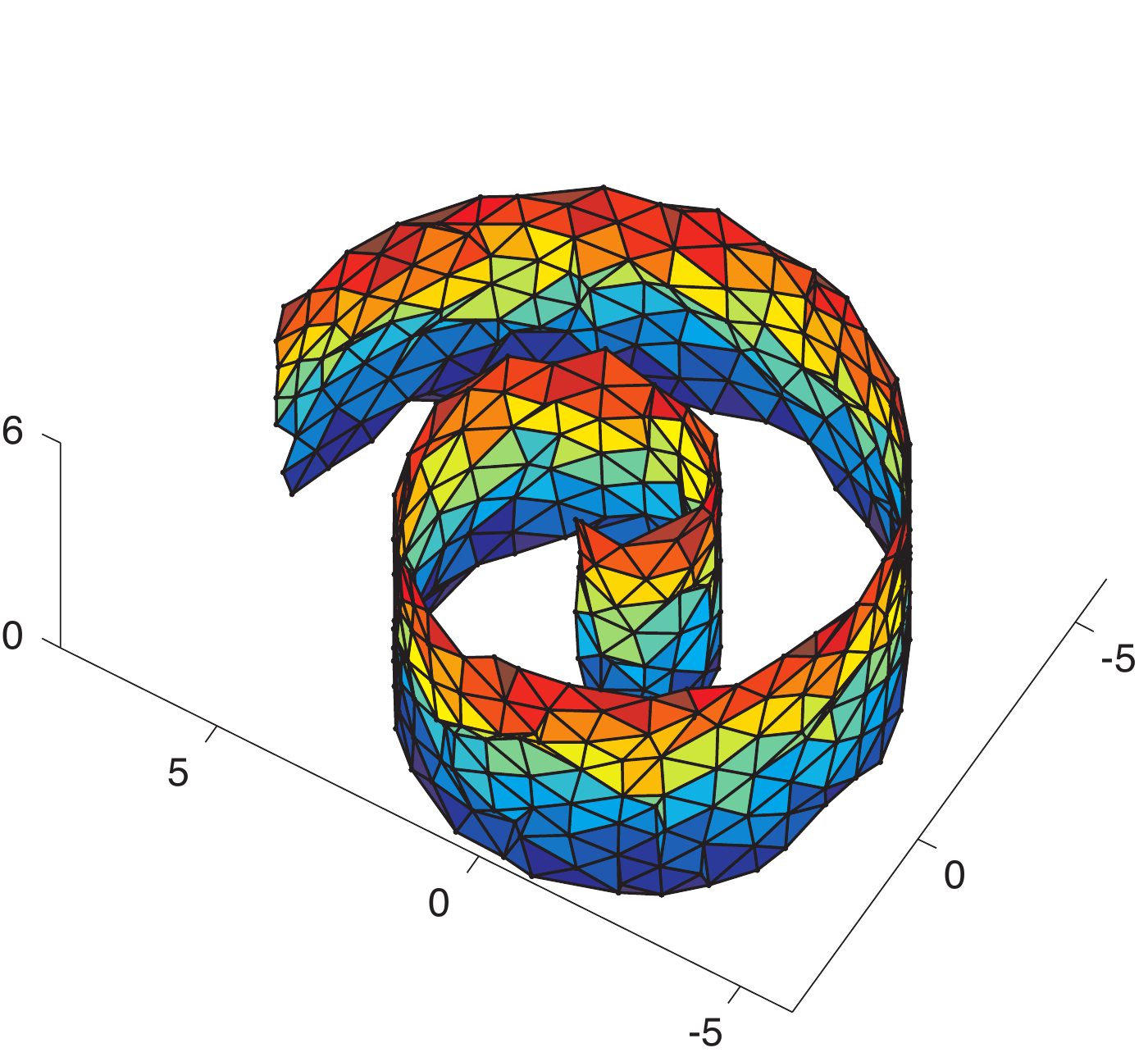}}
%\subfloat[PCA]{}%\includegraphics[scale=.40]{swiss_roll_PCA_gray.pdf}}
\label{fig:swiss_roll}
\caption{We parameterized the swiss roll by $(\tau, \kappa \theta \cos(\theta), \kappa \theta \sin(\theta))$, where $0 \le \tau \le 6$ and $0 \le \theta \le 4 \pi$, and set the curvature parameter $\kappa = 1/2$.
The \emph{notches} on the boundary of the final triangulation account for the large value of $\max_i d_i$.}
\end{center}
\end{figure}

%\clearpage

%%%%%%%%%%%%%%%%%%%%%%%%%%%%%%%%%%%%%%%%%%%%%%%%%%%%%%%%%%%%%%%%%%%%%%%%%%%%%%%%%%%%%%%%%%%%%%%%%%%%%%%%%
%\subsection{Creased sheet}

\begin{figure}[h!]
\begin{center}
\complextable{$1289$}{$669$}{$0.35$}\\
\errtable{$0.10958$}{$0.00083106$}{$0.0055044$}{$10000$}
\subfloat[Data]{\includegraphics[scale=.35]{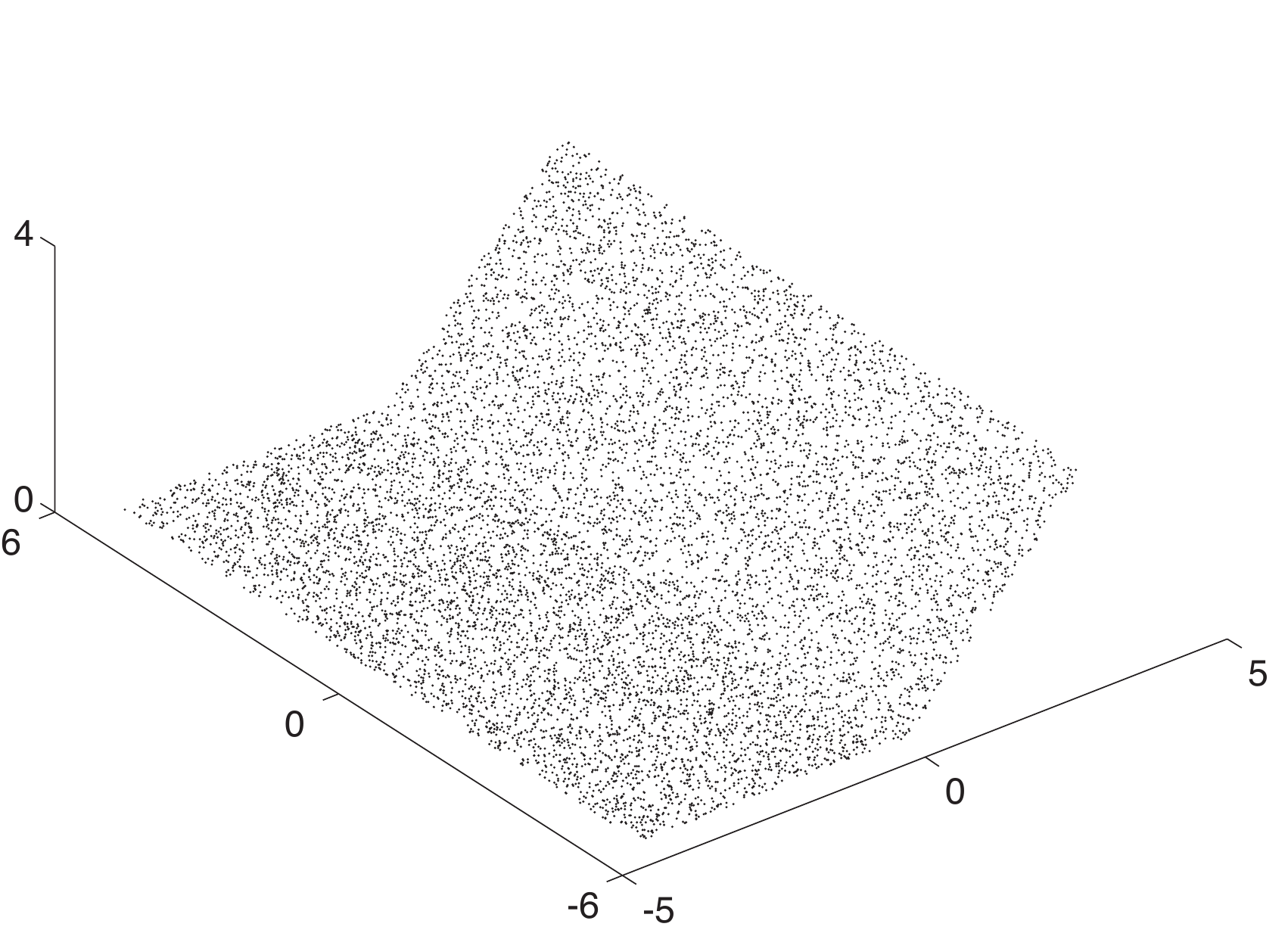}}
\subfloat[SNPCA seams]{\includegraphics[scale=.35]{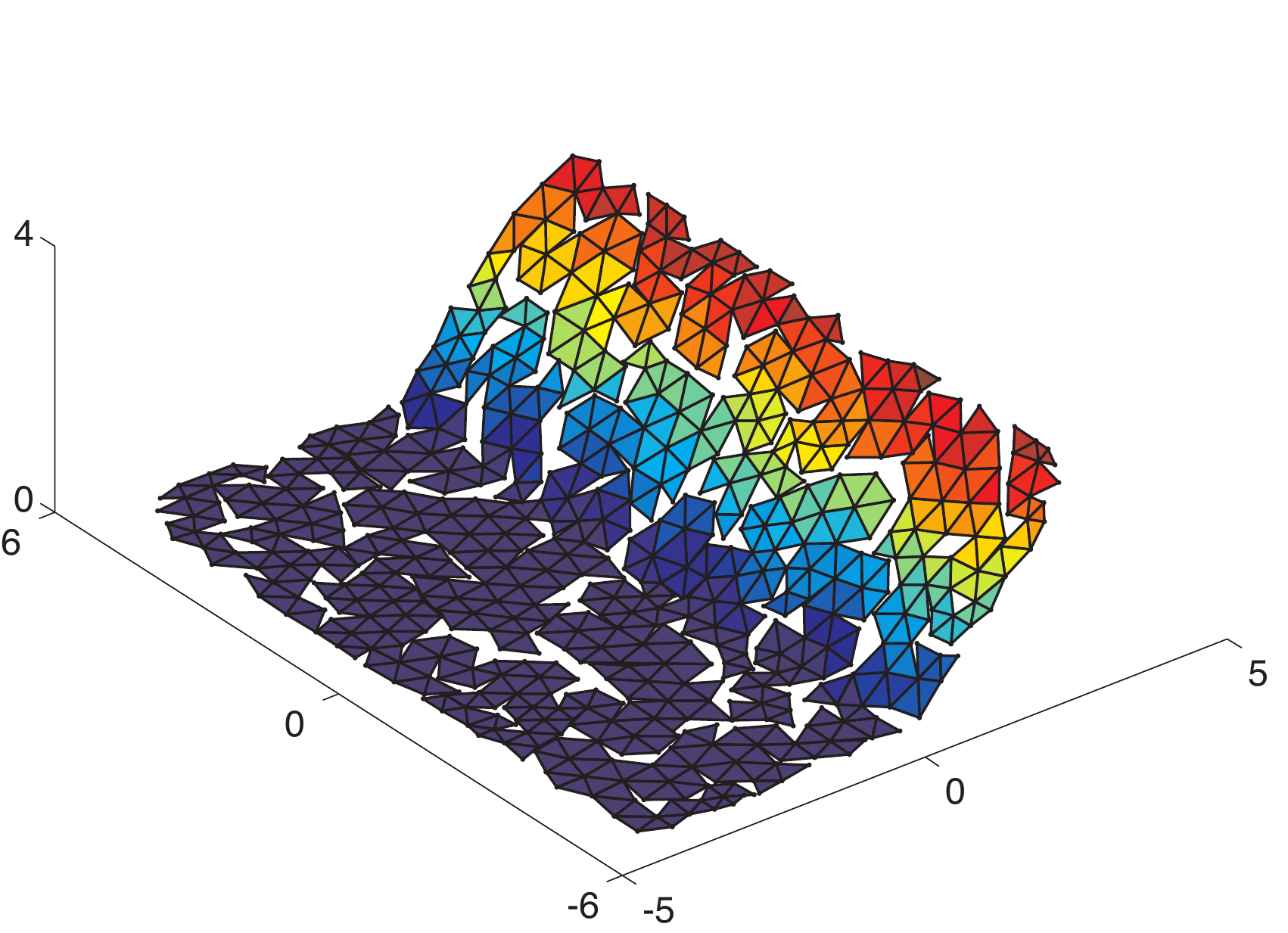}}
\\
\subfloat[SNPCA sewn]{\includegraphics[scale=.35]{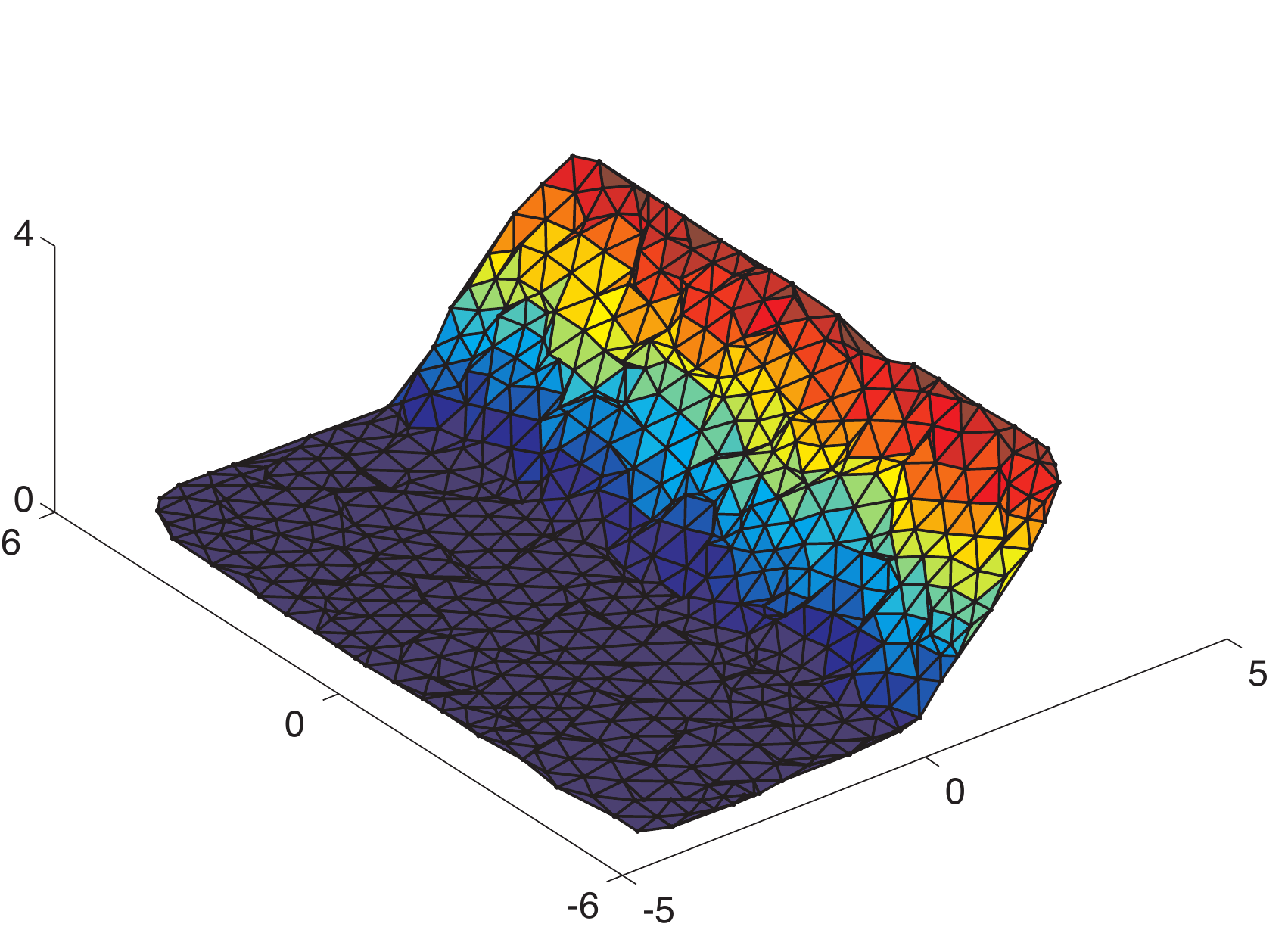}}
%\subfloat[PCA]{}%\includegraphics[scale=.35]{creased_sheet_PCA_gray.pdf}}
\label{fig:creased_sheet}
\caption{The creased sheet data fall on a surface that resembles a creased sheet of $8.5 \times 11$ paper, where the crease angle $.8$ radians.
}
\end{center}
\end{figure}

\clearpage

%%%%%%%%%%%%%%%%%%%%%%%%%%%%%%%%%%%%%%%%%%%%%%%%%%%%%%%%%%%%%%%%%%%%%%%%%%%%%%%%%%%%%%%%%%%%%%%%%%%%%%%%%

%%%%%%%%%%%%%%%%%%%%%%%%%%%%%%%%%%%%%%%%%%%%%%%%%%%%%%%%%%%%%%%%%%%%%%%%%%%%%%%%%%%%%%%%%%%%%%%%%%%%%%%%%
\bibliographystyle{plain}
\bibliography{SNPCA_references}
\end{document}